\documentclass[11pt,a4paper]{article}
\usepackage{multirow}
\usepackage{graphicx}
\usepackage{subfigure}
 \usepackage{epsfig}
 \usepackage{epstopdf}
\usepackage[T1]{fontenc}
\usepackage{extarrows}
\usepackage{geometry}
\usepackage{amsbsy,amsmath,latexsym,amsfonts, epsfig, authblk, amssymb, graphics, bm}
\usepackage{epsf,slidesec,epic,eepic}
\usepackage{fancybox}
\usepackage{fancyhdr}
\usepackage{setspace}
\usepackage{nccmath}
\usepackage{extarrows}
\usepackage{multirow}
\usepackage{extarrows}
\usepackage{graphicx}
\usepackage{algorithm}
\usepackage{algorithmic}
\usepackage{setspace}
\usepackage{subfigure}
\usepackage{graphicx}
\usepackage{amsmath}
\usepackage{booktabs}
\usepackage{multirow}
\usepackage[hidelinks]{hyperref}

\newtheorem{theorem}{Theorem}[section]
\newtheorem{lemma}{Lemma}[section]
\newtheorem{definition}{Definition}[section]

\newtheorem{proposition}{Proposition}[section]
\newtheorem{corollary}{Corollary}[section]
\newtheorem{assumption}{Assumption}
\newtheorem{remark}{Remark}[section]

\usepackage[figuresright]{rotating}
\usepackage{adjustbox}
\usepackage{pgfplots}
\pgfplotsset{compat=1.18}
\usepgfplotslibrary{groupplots}
\newenvironment{proof}[1][Proof]{{\noindent \bf #1:}}{\hfill$\Box$\medskip}

\title{On Primal Superlinear Convergence of Stabilized Sequential Quadratic Programming\thanks{The counterexample in Section~\ref{sec:counterexample} was found with the assistance of GPT-5.6 Sol. AI was also used to polish the language of the manuscript. The authors take full responsibility for the accuracy of this paper.}}

\author{
 Yuqia Wu\footnote{School of Mathematical Sciences, Shenzhen University, 
 	China (\href{mailto:mayuqiawu@szu.edu.cn}{mayuqiawu@szu.edu.cn}).} ,\ \ 
 Yue Wang\footnote{Independent Researcher, ORCID:0000-0001-7988-8791(\href{mailto:yuewangm@outlook.com}{yuewangm@outlook.com}).},\ \
 Yaohua Hu\footnote{School of Mathematical Sciences, Shenzhen University, 
 	China (\href{mailto:mayhhu@szu.edu.cn}{mayhhu@szu.edu.cn}).}
 }
\date{}

\begin{document}
\maketitle
\begin{abstract}
We study primal $Q$-superlinear convergence of stabilized sequential quadratic programming (sSQP). We construct a polynomial equality-constrained problem satisfying the linear independence constraint qualification (LICQ) and the second-order sufficient condition (SOSC), together with an sSQP trajectory whose full primal-dual sequence converges $Q$-quadratically while the primal error increases infinitely often. This gives a negative answer to an open question posed by D.~Fern\'andez, A.~F.~Izmailov, and M.~V.~Solodov [\emph{SIAM J. Optim.}, 20 (2010), pp.~3312--3334]. Under local error bounds and suitable perturbation conditions, we characterize primal $Q$-superlinear convergence of convergent perturbed SQP sequences. For ordinary sSQP, the criterion involves only the primal and multiplier steps.  Based on the criterion, we propose a corrected sSQP algorithm that retains the standard primal step and adaptively corrects the multipliers using only vector operations and scalar comparisons. Under SOSC and locally Lipschitz second derivatives, we establish local $Q$-quadratic convergence of both the primal and full primal-dual sequences, without a constraint qualification or strict complementarity. For equality-constrained problems, SOSC can be replaced by noncriticality of the reference multiplier.
\end{abstract}

\section{Introduction}\label{sec:intro}
Consider the nonlinear programming problem
\begin{equation}\label{eq:nlp}
 \min_{x\in\mathbb{R}^n} f(x)\quad\text{subject to}\quad h(x)=0,\quad g(x)\leq0,
\end{equation}
where $f:\mathbb{R}^n\to\mathbb{R}$, $h:\mathbb{R}^n\to\mathbb{R}^m$ and $g:\mathbb{R}^n\to\mathbb{R}^p$ are twice continuously differentiable near the reference point $\bar{x}$ of interest. The Lagrangian is
$$
 L(x,\lambda,\mu):=f(x)+\lambda^\top h(x)+\mu^\top g(x).
$$
For $z:=(x,\lambda,\mu)$, define the Karush--Kuhn--Tucker (KKT) residual and its norm by
\begin{equation}\label{eq:kkt}
 \mathcal R(z):=
 \begin{pmatrix}
  \nabla_xL(z)\\
  h(x)\\
  \min\{\mu,-g(x)\}
 \end{pmatrix},
 \qquad
 \rho(z):=\|\mathcal R(z)\|.
\end{equation}
The minimum is taken componentwise, so $\rho(z) = 0$ yields $0\le\mu \perp -g(x)\ge0$, where $\perp$ denotes orthogonality. The point $z$ satisfies the KKT conditions if and only if $\rho(z)=0$. We fix a reference KKT pair $\bar z:=(\bar x,\bar\lambda,\bar\mu)$.

At an iterate $z_k:=(x_k,\lambda_k,\mu_k)$, let $H_k:=\nabla^2_{xx}L(z_k)$ and define the quadratic model
$$
 m_k(d):=\nabla f(x_k)^\top d+\tfrac12 d^\top H_kd.
$$
The sequential quadratic programming (SQP) \cite{nocedal2006numerical} subproblem is
\begin{equation}\label{eq:sqp-subproblem}
 \min_{d\in\mathbb R^n} m_k(d)
 \quad\text{subject to}\quad
 h(x_k)+h'(x_k)d=0,\qquad g(x_k)+g'(x_k)d\le0.
\end{equation}
The linearized constraints can be inconsistent, and dependent constraint gradients can lead to nonunique multipliers. Stabilized SQP (sSQP) introduces multiplier variables into the subproblem~\cite{wright1998superlinear}. Let $\rho_k:=\rho(z_k)=\|\mathcal R(z_k)\|$. For $\rho_k>0$, the sSQP subproblem is \begin{equation}\label{eq:ssqp-subproblem}
\begin{aligned}
 \min_{d,\lambda,\mu}\quad
 &m_k(d)+\tfrac{\rho_k}{2}
       \bigl(\|\lambda\|^2+\|\mu\|^2\bigr),\\
 \text{subject to}\quad
 &h(x_k)+h'(x_k)d-\rho_k(\lambda-\lambda_k)=0,\\
 &g(x_k)+g'(x_k)d-\rho_k(\mu-\mu_k)\le0.
\end{aligned}
\end{equation}
The additional multiplier variables make the constraints of~\eqref{eq:ssqp-subproblem} feasible for every $\rho_k>0$~\cite{fernandez2010stabilized}. Given a stationary point $(d_k,\lambda_{k+1},\mu_{k+1})$ of~\eqref{eq:ssqp-subproblem}, the next primal iterate is $x_{k+1}=x_k+d_k$.

Wright~\cite{wright1998superlinear} established local superlinear convergence of the primal-dual sSQP sequence under the Mangasarian--Fromovitz constraint qualification (MFCQ) and a second-order sufficient condition (SOSC) holding uniformly over the multiplier set, with the initial multiplier sufficiently close to a strictly complementary multiplier. Hager~\cite{hager1999stabilized} removed the MFCQ and strict complementarity requirements by assuming the strong second-order sufficient condition at a reference primal-dual solution. Fern\'andez and Solodov~\cite{fernandez2010stabilized} subsequently weakened this condition to the usual SOSC. Under SOSC, the full primal-dual sequence converges locally $Q$-superlinearly without any constraint qualification or strict complementarity, and the rate is strengthened to $Q$-quadratic when the second-order derivatives are locally Lipschitz continuous. For equality-constrained problems, Izmailov and Solodov~\cite{izmailov2012stabilized} further showed that SOSC can be replaced by noncriticality of the reference multiplier.

The development of sSQP has also focused on Hessian approximation and globalization techniques. Quasi-Newton updates were studied in~\cite{fernandez2013quasi}. Globalization strategies have been developed using primal-dual augmented Lagrangian functions and line searches~\cite{gill2013globally,gill2017stabilizeda,gill2017stabilizedb}, augmented Lagrangian steps~\cite{izmailov2015combining}, and smooth exact penalty functions~\cite{izmailov2016globalizing}. Stabilized methods have also been extended to nonlinear semidefinite programming~\cite{yamakawa2022stabilized}, Banach-space optimization~\cite{yamakawa2026local}, and composite optimization~\cite{hang2026stabilized}. For a review of stabilized methods, see~\cite{fernandez2014stabilized,izmailov2014newton}.

In optimization, the convergence rate of the primal sequence is of particular interest. For a sequence $z_k\to\bar z$, primal-dual $Q$-superlinear convergence, $\|z_{k+1}-\bar z\|=o(\|z_k-\bar z\|)$, does not in general imply primal $Q$-superlinear convergence, $\|x_{k+1}-\bar x\|=o(\|x_k-\bar x\|)$. Primal convergence has been established for several constrained optimization methods. Fern\'andez, Izmailov and Solodov~\cite{fernandez2010sharp} obtained sharp primal superlinear convergence results in a general perturbed SQP framework. For ordinary SQP, Izmailov and Solodov~\cite[Theorem~2.2]{izmailov2016new} proved a primal superlinear estimate from a sufficiently close primal-dual initialization, under multiplier uniqueness and SOSC. For the augmented Lagrangian method, Fern\'andez and Solodov~\cite{fernandez2012local} showed that, $Q$-superlinear convergence of the primal-dual sequence implies $Q$-superlinear convergence of the primal sequence under SOSC. Primal superlinear convergence of SQP has also been studied in composite optimization by Sarabi~\cite{sarabi2022primal}.

For sSQP, however, primal $Q$-superlinear convergence under the usual local convergence assumptions was left as an open question in~\cite[Section~4.4]{fernandez2010sharp}. The sSQP iteration can be regarded as a perturbed SQP (pSQP) step in the framework of~\cite{izmailov2010inexact,fernandez2010sharp}. Write
$$
 d_k:=x_{k+1}-x_k,\qquad
 \Delta\lambda_k:=\lambda_{k+1}-\lambda_k,\qquad
 \Delta\mu_k:=\mu_{k+1}-\mu_k.
$$
The subproblem of the pSQP reads
\begin{equation}\label{eq:psqp}
\begin{aligned}
0={}&\nabla_xL(z_k)+H_kd_k
       +h'(x_k)^\top\Delta\lambda_k
       +g'(x_k)^\top\Delta\mu_k+\omega_1^k,\\
0={}&h(x_k)+h'(x_k)d_k+\omega_2^k,\\
0\le{}&\mu_{k+1}\ \perp\
       -g(x_k)-g'(x_k)d_k-\omega_3^k\ge0,
\end{aligned}
\end{equation}
where $\omega_1^k\in\mathbb R^n$, $\omega_2^k\in\mathbb R^m$, and $\omega_3^k\in\mathbb R^p$ are perturbations. Setting all three perturbations to zero recovers the KKT conditions of~\eqref{eq:sqp-subproblem}. The ordinary sSQP equations correspond to
\begin{equation}\label{eq:ssqp-omega}
 \omega_1^k=0,\qquad
 \omega_2^k=-\rho_k\Delta\lambda_k,\qquad
 \omega_3^k=-\rho_k\Delta\mu_k.
\end{equation}

For ordinary sSQP, the perturbations~\eqref{eq:ssqp-omega} involve the multiplier increments. The usual local convergence assumptions do not ensure that these perturbations are sufficiently small relative to the primal step to apply the results in~\cite{fernandez2010sharp}. The example in Section~4.4 of that paper shows that a single sSQP step can increase the primal error even under LICQ and SOSC, but the subsequent iterates recover primal $Q$-superlinear convergence. Thus a single increase does not settle the question. A related result of Fern\'andez \cite[Proposition~5]{fernandez2013quasi} establishes the two-step estimate:
\begin{equation}\label{eq:intro-two-step}
 \frac{\|x_{k+1}-\bar x\|}
 {\max\{\|x_k-\bar x\|,\|x_{k-1}-\bar x\|\}}
 \rightarrow0.
\end{equation}
This estimate allows individual primal errors to increase, so it does not imply the one-step rate.

\subsection{Contributions}
In this paper, we first give a negative answer to the open question raised in~\cite[Section~4.4]{fernandez2010sharp}. We consider a polynomial equality-constrained problem satisfying LICQ and SOSC, and construct an sSQP trajectory whose full primal-dual sequence converges $Q$-quadratically while the primal error increases for infinite times. This sequence also satisfies~\eqref{eq:intro-two-step}, showing that the two-step estimate cannot in general be strengthened to one-step primal $Q$-superlinear convergence.

Under local error bounds and suitable perturbation conditions, we characterize the sufficient and necessary conditions for primal $Q$-superlinear convergence of the pSQP sequences. For ordinary sSQP under a local primal error bound, the criterion reduces to $\|(\Delta\lambda_k,\Delta\mu_k)\|^2=o(s_k)$, where $s_k:=\|x_{k+1}-x_k\|$.  Moreover, with locally Lipschitz second derivatives, we prove that the primal $Q$-quadratic convergence is equivalent to $\rho_k=O(s_k)$.

The criterion motivates an adaptive multiplier correction. After solving the ordinary sSQP subproblem, we retain its primal point and use the residual norm $\rho_k$ and the primal step length $s_k$ to choose the weight in the multiplier update. Under SOSC and locally Lipschitz second derivatives, we prove local $Q$-quadratic convergence of both the primal and primal-dual sequences, without a constraint qualification or strict complementarity. For equality-constrained problems, the same conclusions hold when SOSC is replaced by noncriticality of the reference multiplier.

Finally, numerical experiments on the counterexample illustrate the failure of one-step primal superlinear convergence and the effect of the multiplier correction. We also compare ordinary and corrected sSQP on selected problems from the DEGEN collection in terms of local convergence behavior and computational cost.

\medskip
\noindent {\bf Notations.} Throughout this paper, $\mathbb{R}^n$ denotes the $n$-dimensional Euclidean space equipped with standard inner product $\langle \cdot, \cdot \rangle$ and its induced norm $\|\cdot\|,$ and matrix norms are the corresponding operator norms.  We denote the closed Euclidean ball centered at $u$ with radius $r>0$ by $\mathbb B(u,r):=\{v:\|v-u\|\le r\}$. For a nonempty set $S$, the distance from $u$ to $S$ is $\operatorname{dist}(u,S):=\inf_{v\in S}\|u-v\|$. Subscripts indexed by a set select the corresponding vector components or matrix rows.  Vector inequalities, minima, and maxima are interpreted componentwise, and $[a]_+:=\max\{a,0\}$ denotes the positive part of $a$. The complementarity notation $0\le a\perp b\ge0$ means that $a\ge0$, $b\ge0$, and $a^\top b=0$.

\section{Preliminaries}\label{sec:prelim}

We recall the local error bounds used in the convergence analysis and the conditions that ensure local solvability of the sSQP subproblem. We also relate the size of an sSQP step to the KKT residual at its input.
\subsection{KKT residual and error bound}

Recall that $\bar z=(\bar x,\bar\lambda,\bar\mu)$ is a reference KKT point. For $x$ near $\bar x$, write
$$
 A(x):=h'(x),\quad G(x):=g'(x),\quad
 \bar A:=A(\bar x),\quad \bar G:=G(\bar x),\quad
 \bar H:=\nabla^2_{xx}L(\bar z).
$$
We define
\begin{equation*}
\begin{aligned}
  \mathcal A:=\{i:g_i(\bar x)=0\},&\quad \mathcal N:=\{1,\ldots,p\}\setminus\mathcal A, \\
  \mathcal A_+:=\{i\in\mathcal A:\bar\mu_i>0\},& \quad \mathcal A_0:=\mathcal A\setminus\mathcal A_+.
 \end{aligned}
 \end{equation*}
We say that the linear independence constraint qualification (LICQ) holds at $\bar x$ if the vectors
$$
 \nabla h_i(\bar x),\  i=1,\ldots,m, \quad
 \nabla g_i(\bar x), \ i\in\mathcal A
$$
are linearly independent. LICQ guarantees uniqueness of the multiplier at $\bar x$.  The curvature condition used in the local convergence analysis of sSQP is imposed on the critical cone
\begin{equation}\label{eq:critical-cone}
 \mathcal{C}:=\{v:\bar A v=0,\ \bar G_{\mathcal A_+}v=0,
                 \ \bar G_{\mathcal A_0}v\le0\}.
\end{equation}
The critical cone consists of the directions that satisfy the linearized constraints and along which the objective has no first-order change. Along these directions, first-order information alone does not determine local optimality, so second-order information is needed.  The second-order sufficient condition (SOSC) requires positive curvature of the Lagrangian on every nonzero direction in $\mathcal{C}$.

\begin{definition}
Problem~\eqref{eq:nlp} satisfies the second-order sufficient condition (SOSC) at $\bar z$ if
\begin{equation*}
 v^\top\bar H v>0\qquad\text{for every }0\ne v\in \mathcal{C}.
\end{equation*}
\end{definition}

For problems with nonunique multipliers, the convergence rate of Newton-type methods may depend on the multiplier approached by the dual sequence. Attraction to critical multipliers and the resulting loss of convergence speed were studied in~\cite{izmailov2009attraction,izmailov2009examples,izmailov2011attraction}; see also the surveys~\cite{izmailov2015critical,fernandez2014stabilized}. We recall the definitions of critical and noncritical multipliers.

\begin{definition}[Critical and noncritical multipliers]
The multiplier $(\bar\lambda,\bar\mu)$ is called critical if there exist $v\in\mathbb R^n\setminus\{0\}$, $w\in\mathbb R^m$, and $\eta\in\mathbb R^p$ satisfying
\begin{equation}\label{eq:noncritical-system}
\begin{aligned}
 \bar H v+\bar A^\top w+\bar G^\top\eta&=0,\qquad \bar A v=0,\\
 \bar G_{\mathcal A_+}v&=0,\qquad \eta_{\mathcal N}=0,\\
 0\le\eta_{\mathcal A_0}&\ \perp\ -\bar G_{\mathcal A_0}v\ge0.
\end{aligned}
\end{equation}
Otherwise, the multiplier is called noncritical.
\end{definition}

As shown in Fern\'andez, Izmailov and Solodov~\cite{fernandez2010sharp}, SOSC implies noncriticality. Indeed, taking the inner product of the stationarity equation in~\eqref{eq:noncritical-system} with $v$ gives $v^\top\bar H v=0$. Since $v\in \mathcal{C}$, SOSC then yields $v=0$.  The relaxation from SOSC to noncriticality is useful in the local analysis of sSQP. For equality-constrained problems, Izmailov and Solodov~\cite{izmailov2012stabilized} established local superlinear convergence of the primal-dual sSQP sequence under noncriticality.

The connection between noncriticality and error bounds also plays a role in the analysis of primal superlinear convergence in~\cite{fernandez2010sharp}. Related results on perturbation stability are developed in~\cite{izmailov2013note}. For problems with both equality and inequality constraints, write $y=(\lambda,\mu)$ and let
\begin{equation}\label{eq:defM}
 \mathcal M_*:=\{(\lambda,\mu):\nabla_xL(\bar x,\lambda,\mu)=0,
       \ \mu_{\mathcal A}\ge0,\ \mu_{\mathcal N}=0\}
\end{equation}
be the multiplier set associated with $\bar x$. Izmailov and Solodov~\cite[Proposition~3 and Remark~2]{izmailov2012stabilized} show that, $(\bar\lambda,\bar\mu)$ is noncritical if and only if there exists $K>0$ such that
\begin{equation}\label{eq:primal-dual-eb}
 \|x-\bar x\|+\operatorname{dist}(y,\mathcal M_*)
 \le K\|\mathcal R(z)\|
 \qquad\text{for $z$ sufficiently close to $\bar z$}.
\end{equation}
Such estimates are also used in the general local convergence frameworks for nonisolated solutions and perturbed Newton steps developed in~\cite{fischer2002local,izmailov2010inexact}. The primal part of~\eqref{eq:primal-dual-eb} allows us to estimate the error at the next iterate from its KKT residual. We formulate the primal error bound under the following assumption.
\begin{assumption}[Primal error bound]\label{ass:eb}
There exist a neighborhood $U$ of $\bar z$ and $\kappa>0$ such that
\begin{equation}\label{eq:ordinary-eb}
 \|x-\bar x\|\le\kappa\|\mathcal R(z)\|
 \qquad\text{for every $z\in U$}.
\end{equation}
\end{assumption}

Assumption~\ref{ass:eb} therefore holds under SOSC or, more generally, under noncriticality. This condition is essential in characterizing the primal rate of the sSQP sequence.  We further assume Lipschitz continuity of the second derivatives of $f, g, h$.

\begin{assumption}\label{ass:lipschitz}
There exist $\delta_0>0$ and constants $L_\varphi>0$ such that, for all $x,x'\in\mathbb B(\bar x,\delta_0)$,
$$
 \|\nabla^2\varphi(x)-\nabla^2\varphi(x')\|
 \le L_\varphi\|x-x'\|,
 \qquad \varphi\in\{f,h_1,\ldots,h_m,g_1,\ldots,g_p\}.
$$
\end{assumption}

Under Assumption~\ref{ass:lipschitz}, we may reduce $\delta_0$ and choose $M\ge1$ such that, for all $x,x'\in\mathbb B(\bar x,\delta_0)$ and $y,y'\in\mathbb B((\bar\lambda,\bar\mu),\delta_0)$,
\begin{equation}\label{eq:local-smoothness-bounds}
\begin{aligned}
 \|A(x)\|+\|G(x)\|+\|\nabla^2_{xx}L(x,y)\|
 &\le M,\\
 \|A(x)-A(x')\|+\|G(x)-G(x')\|
 &\le M\|x-x'\|,\\
 \|\nabla^2_{xx}L(x,y)-\nabla^2_{xx}L(x',y)\|
 &\le M\|x-x'\|,\\
 \|\mathcal R(x,y)-\mathcal R(x',y')\|
 &\le M\|(x,y)-(x',y')\|.
\end{aligned}
\end{equation}

\subsection{Local sSQP steps}
For the local estimates below, write $z=(x,\lambda,\mu)$ and $z^+=(x^+,\lambda^+,\mu^+)$ for the input and output of one sSQP step. Set
\begin{equation}\label{eq:notions}
\begin{aligned}
 A:=A(x),\quad G:=G(x),\quad H:=\nabla^2_{xx}L(z),\quad
 \rho:=\rho(z)>0,\\
 d:=x^+-x,\qquad
 \Delta\lambda:=\lambda^+-\lambda,\qquad
 \Delta\mu:=\mu^+-\mu.
\end{aligned}
\end{equation}
In this notation, the KKT system of sSQP subproblem \eqref{eq:ssqp-subproblem} is
\begin{equation}\label{eq:ssqp-system}
\begin{aligned}
0={}&\nabla_xL(z)+Hd+A^\top\Delta\lambda+G^\top\Delta\mu,\\
0={}&h(x)+Ad-\rho\Delta\lambda,\\
0\leq{}&\mu^+\ \perp\ -g(x)-Gd+\rho\Delta\mu\geq0.
\end{aligned}
\end{equation}

Local solvability is a key ingredient in the superlinear convergence analysis of sSQP. The standard local analysis combines an error bound for the KKT system with the existence of nearby subproblem stationary points whose primal-dual steps are controlled by the current residual. When inequality constraints are present, noncriticality alone may fail to ensure local solvability and therefore cannot, in general, replace SOSC in the superlinear convergence theory~\cite{izmailov2012stabilized}. The following proposition gives a bound on the primal-dual step for nearby solutions of the sSQP system and states conditions ensuring the existence of such solutions.

\begin{proposition}[Local step bound and solvability]\label{prop:local-solvability}
The following assertions hold.
\begin{itemize}
\item[{\rm (i)}] Under Assumption~\ref{ass:eb}, there exist a neighborhood $U$ of $\bar z$ and $C_1>0$ such that every pair $z,z^+\in U$ satisfying the sSQP system~\eqref{eq:ssqp-system}, with $\mu\geq0$ and $\rho>0$, obeys
\begin{equation}\label{eq:zbound}
 \|z^+-z\|\leq C_1\rho.
\end{equation}
\item[{\rm (ii)}] Suppose that either SOSC holds at $\bar z$, or the problem has only equality constraints and $\bar\lambda$ is noncritical. Let $U$ and $C_1$ be as in part~{\rm (i)}. For every input $z$ sufficiently close to $\bar z$ with $\mu\geq0$ and $\rho>0$, the sSQP system~\eqref{eq:ssqp-system} admits a solution $z^+\in U$ satisfying \eqref{eq:zbound}. In the equality-constrained case, this solution is unique.
\end{itemize}
\end{proposition}
\begin{proof} \ 
Part~{\rm (i)} is proved in Appendix~\ref{app:small-step}. For part~{\rm (ii)}, local solvability under SOSC follows from~\cite[Theorem~3]{fernandez2010stabilized}. For equality constraints, local solvability and uniqueness follow under noncriticality by~\cite[Lemma~2 and Proposition~2]{izmailov2012stabilized}, which yields $\|z^+-z\|=O(\rho)$. After shrinking the input neighborhood and enlarging $C_1$ if necessary, the selected solution lies in $U$ and satisfies~\eqref{eq:zbound}.
\end{proof}

We remark here that part~{\rm (i)} provides a uniform step bound for every nearby solution of the sSQP system, without asserting its existence. Part~{\rm (ii)} guarantees local existence under additional assumptions, as well as uniqueness in the equality-constrained case. For the reduced sSQP system, a step bound of the form~\eqref{eq:zbound} was established under SOSC in~\cite[Proposition~3]{fernandez2010stabilized}. Here, part~{\rm (i)} establishes this bound for all nearby pairs $z,z^+$ satisfying the full sSQP system under the primal error bound in Assumption~\ref{ass:eb}, which is weaker than SOSC.

Next, we consider the reverse estimate of equation \eqref{eq:zbound}, which  does not require an error bound condition.

\begin{lemma}\label{lem:upper}
There exist a neighborhood $U$ of $\bar z$ and $C_2>0$ such that every pair $z,z^+\in U$ satisfying the sSQP system~\eqref{eq:ssqp-system} obeys $\rho\leq C_2\|z^+-z\|.$
\end{lemma}
\begin{proof}\
Choose a bounded neighborhood $U$ on which $H$, $A$, $G$ and $\rho$ defined in \eqref{eq:notions} are uniformly bounded. Fix $C>0$ such that $\|H\|+2\|A\|+\|G\|+\rho+1\le C$ throughout $U$. The first two equations in \eqref{eq:ssqp-system} give
$$
 \|\nabla_xL(z)\|+\|h(x)\|\le C\|z^+-z\|.
$$
The complementarity equation in \eqref{eq:ssqp-system} along with the Lipschitz continuity of minimum function yields
$$
 \|\min\{\mu,-g(x)\}\|
 \le\|\Delta\mu\|+\|Gd-\rho\Delta\mu\|
 \le C\|z^+-z\|.
$$
Adding the two inequalities proves the result with $C_2=2C$.
\end{proof}

Under Assumption~\ref{ass:eb}, Proposition~\ref{prop:local-solvability}{\rm (i)} and Lemma~\ref{lem:upper} show that
$$
 C_2^{-1}\rho\le\|z^+-z\|\le C_1\rho
$$
uniformly for all $z,z^+$ sufficiently close to $\bar z$ that satisfy the sSQP system~\eqref{eq:ssqp-system} with $\mu\ge0$ and $\rho>0$.

\section{Failure of primal $Q$-superlinear convergence}\label{sec:counterexample}

In this section, we construct a counterexample to one-step primal $Q$-superlinear convergence of sSQP. The example satisfies LICQ and SOSC and has a single equality constraint. Along the constructed trajectory, the primal error increases infinitely often, while the full primal-dual sequence converges $Q$-quadratically.

Consider the following problem:
\begin{equation}\label{eq:example}
\begin{aligned}
  \min_{(x,y)\in\mathbb{R}^2}\quad & F(x,y):=\frac{x^3}{3}+\frac{y^2}{2},\\
  \text{subject to}\quad &h(x,y):=x=0.
 \end{aligned}
\end{equation}
Let $\lambda \in \mathbb{R}$ be the Lagrangian multiplier. The Lagrangian function of \eqref{eq:example} is
$$ L(x, y; \lambda) = \frac{x^3}{3}+\frac{y^2}{2} + \lambda x,$$
and its derivatives with respect to the primal variables are
\begin{equation*}
 \nabla L(x,y)=(x^2+\lambda,y)^\top,\qquad
 \nabla^2 L(x,y)=\begin{bmatrix}2x&0\\0&1\end{bmatrix},\qquad
 \nabla h(x,y)=(1,0)^\top.
\end{equation*}

Feasibility gives $x=0$, so the unique global solution is $(0,0)$ and its unique multiplier is zero. Since $h'(0,0)=(1,0)$, LICQ holds. On the tangent space $\{(0,d_y):d_y\in\mathbb{R}\}$, the Lagrangian Hessian has quadratic form $d_y^2$, so SOSC holds.

Next we describe the sSQP iterates. The KKT system for problem~\eqref{eq:example} is
$$
 x^2+\lambda=0,\quad y=0,\quad x=0,
$$
and the corresponding residual function at $(x,y,\lambda)$ is
\begin{equation*}
 \rho(x,y,\lambda)=\sqrt{(x^2+\lambda)^2+y^2+x^2}.
\end{equation*}
The classical sSQP subproblem solves the following to obtain $(x^+, y^+, \lambda^+)$:
\begin{equation}\label{eq:ssqp-qp}
\begin{aligned}
 \min_{(x^+, y^+, \lambda^+) \in \mathbb{R}^3}\quad
 &(x^2,y)\binom{x^+-x}{y^+-y}
 +x(x^+-x)^2+\frac12(y^+-y)^2+\frac{\rho(x,y,\lambda)}{2}(\lambda^+)^2,\\
 \text{subject to}\quad & x^+ -\rho(x,y,\lambda)(\lambda^+-\lambda)=0,
\end{aligned}
\end{equation}
The stationarity equation for $y^+$ gives $y^+=0$, so every trajectory lies in the subspace $y=0$ after its first step. On this subspace, we identify $(x,0,\lambda)$ with $(x,\lambda)$ and write $\rho:=\rho(x,\lambda):=\sqrt{(x^2+\lambda)^2+x^2}$. For every nonzero $(x,\lambda)$ sufficiently close to the origin, we have $\rho>0$ and $1+2x\rho>0$. Solving the remaining KKT equations of~\eqref{eq:ssqp-qp} gives the reduced iteration map
\begin{equation}\label{eq:exact-map}
 \Phi(x,\lambda):=(x^+,\lambda^+), \quad {\rm where}\ \ x^+=\frac{(x^2-\lambda)\rho}{1+2x\rho},
 \qquad
 \lambda^+=x^2-2xx^+.
\end{equation}

The factor $x^2-\lambda$ in~\eqref{eq:exact-map} suggests choosing $\lambda$ close to $x^2$. The first step can then make $|x^+|$ much smaller than $|\lambda^+|$, with $\lambda^+\approx x^2$. In this case, the next step gives $x^{++}\approx-(\lambda^+)^2$, so $|x^{++}|$ can exceed $|x^+|$. To make this increase recur along a single trajectory, we seek points for which the relation $\lambda\approx x^2$ is recovered after two steps. For $t>0$ and $a>0$, consider the primal-dual points

\begin{equation}\label{eq:parapoint}
 z(t,a):=\bigl(-t,t^2(1-ta)\bigr).
\end{equation}
For $(x,\lambda)=z(t,a)$, we have
\begin{equation}\label{eq:xlambda}
 x^2-\lambda=at^3,\qquad
 \rho\bigl(z(t,a)\bigr)
 =t\sqrt{1+t^2(2-ta)^2}.
\end{equation}

The next lemma analyzes two consecutive sSQP iterations starting from the points $z(t,a)$ defined in \eqref{eq:parapoint}, whose proof is given in Appendix~\ref{app:curve}.

\begin{lemma}\label{lem:two-step}
Let $\Phi$ be the reduced sSQP map defined by~\eqref{eq:exact-map}, and $z(t,a)$ be given by~\eqref{eq:parapoint}. Set $c:=\sqrt2-1$ and fix a compact interval $I\subset (0,\infty)$ with $c\in {\rm int}(I)$. Then there exist $\delta_0>0$ and functions $T=T(t,a)$ and $R=R(t,a)$, jointly smooth on an open neighborhood of $[0,\delta_0]\times I$, such that the following equations hold for $0\le t\le\delta_0$ and $a\in I$:
\begin{subequations}
\begin{align}
 \Phi(z(t,a))
 &=\bigl(t^4(a+O(t^2)),t^2(1+O(t^3))\bigr),
 \label{eq:first-step}\\
 \Phi^2(z(t,a))
 &=\bigl(-T(t,a),T(t,a)^2(1+R(t,a))\bigr),
 \label{eq:two-arrows}\\
 T(t,a)
 &=t^4\bigl(1+t^3\theta(t,a)\bigr),
 \label{eq:T-expansion}\\
 R(t,a)
 &=a^2+2a-1+t^2\eta(t,a),
 \label{eq:R-expansion}
\end{align}
\end{subequations}
where $\theta$ and $\eta$ are jointly smooth on the same neighborhood and bounded on $[0,\delta_0]\times I$. In particular,
\begin{equation}\label{eq:T-derivatives}
 \partial_aT(t,a)=O(t^7),\qquad
 \partial_tT(t,a)=4t^3+O(t^6),
\end{equation}
uniformly for $a\in I$.
\end{lemma}

By~\eqref{eq:two-arrows}, two sSQP iterations send the first component $-t$ to $-T(t,a)$, so the curve parameter is updated from $t$ to $T(t,a)$. To construct an invariant curve, write $a=q(t)$ and set $t_+=T(t,q(t))$. For the image to lie on the same curve, its second component must equal $t_+^2(1-t_+q(t_+))$. Comparing this with $t_+^2(1+R(t,q(t)))$ in~\eqref{eq:two-arrows} gives
\begin{equation}\label{eq:invariance}
 R(t,q(t))+T(t,q(t))q(T(t,q(t)))=0.
\end{equation}
At $t=0$, the term involving $T$ vanishes, so the invariance equation requires $R(0,q(0))=0$. By~\eqref{eq:R-expansion},
\begin{equation}\label{eq:R-zero}
 R(0,a)=a^2+2a-1,\qquad
 R(0,c)=0,\qquad
 \partial_aR(0,c)=2\sqrt2>0.
\end{equation}
Using $\partial_aR(0,c)>0$, we construct successive approximations to a solution of~\eqref{eq:invariance} by the implicit function theorem. Their limit defines the invariant curve in the following lemma. The proof is given in Appendix~\ref{app:curve}.

\begin{lemma}\label{lem:curve}
Let $\Phi$, $T$, $R$, and $c$ be as in Lemma~\ref{lem:two-step}. Then there exist $\delta>0$ and a Lipschitz continuous function $q_*:[0,\delta]\to\mathbb R$, satisfying $q_*(0)=c$, such that the following assertions hold:
\begin{itemize}
\item[{\rm (i)}] Define
$$
 \mathcal C_{q_*}
 :=\{z(t,q_*(t)):0<t\leq\delta\}
 =\{(-t,t^2(1-tq_*(t))):0<t\leq\delta\}.
$$
For every $t\in(0,\delta]$, set $t_+:=T(t,q_*(t))$. Then
$$
 0<t_+<t,\quad
 \Phi^2(z(t,q_*(t)))
 =z(t_+,q_*(t_+))\in\mathcal C_{q_*}.
$$

\item[{\rm (ii)}] With $t_+$ defined in {\rm (i)}, it holds that $t_+=t^4(1+o(1)).$
\end{itemize}
\end{lemma}

We now use the invariant curve in Lemma~\ref{lem:curve} to construct an sSQP trajectory whose primal error increases infinitely often.

\begin{theorem}\label{thm:counterexample}
Let $\delta>0$ and $q_*:[0,\delta]\to\mathbb R$ be given by Lemma~\ref{lem:curve}, with $\delta$ reduced if necessary. Choose any $t_0\in(0,\delta]$ and start the sSQP method for \eqref{eq:example} from
\begin{equation*}
 (x_0,y_0,\lambda_0)
 =\bigl(-t_0,0,t_0^2(1-t_0q_*(t_0))\bigr).
\end{equation*}
Then the following assertions hold:
\begin{itemize}
\item[{\rm (i)}] Every sSQP subproblem \eqref{eq:ssqp-qp} has a unique global solution.
\item[{\rm (ii)}] The primal-dual sequence converges $Q$-quadratically to $(0,0,0)$.
\item[{\rm (iii)}] The primal error $e_k:=\|(x_k,y_k)\|$ satisfies $e_{2j+2}/e_{2j+1}\to1+\sqrt{2}>1$. Consequently, the primal sequence fails to converge $Q$-superlinearly.
\item[{\rm (iv)}] The primal sequence satisfies the two-step estimate: $ \frac{e_{k+1}}{\max\{e_k, e_{k-1}\}} \rightarrow0. $
\end{itemize}
\end{theorem}

\begin{proof} \ 
Set $a_j:=q_*(t_j),\ t_{j+1}:=T(t_j,a_j).$ By~\eqref{eq:T-expansion}, reducing $\delta$ if necessary, we have $0<t_{j+1}\le 2t_j^4\le t_j/2$. Thus $t_j\to0$, and Lemma~\ref{lem:curve} gives
\begin{equation}\label{eq:t-recurrence}
 t_{j+1}=t_j^4(1+o(1)),\qquad a_j = q_*(t_j)\rightarrow c.
\end{equation}
Since $y_0=0$ and the sSQP update gives $y^+=0$, we have $y_k=0$ for all $k$. The invariance of $\mathcal C_{q_*}$ under $\Phi^2$, together with~\eqref{eq:first-step}, gives
\begin{equation}\label{eq:iterates}
\begin{aligned}
 (x_{2j},y_{2j},\lambda_{2j})
 &=\bigl(-t_j,0,t_j^2(1-t_ja_j)\bigr), \\
 (x_{2j+1},y_{2j+1},\lambda_{2j+1})
 &=\bigl(t_j^4(a_j+O(t_j^2)),0, t_j^2(1+O(t_j^3))\bigr), \\
 (x_{2j+2},y_{2j+2},\lambda_{2j+2})
 &=\bigl(-t_{j+1},0, t_{j+1}^2(1-t_{j+1}a_{j+1})\bigr).
\end{aligned}
\end{equation}

\noindent {\rm (i)} Eliminating $x^+=\rho(\lambda^+-\lambda)$ from~\eqref{eq:ssqp-qp} gives a quadratic function in $(y^+,\lambda^+)$ with Hessian $\operatorname{diag}(1,(1+2x\rho)\rho)$. For $\delta$ sufficiently small, $\rho>0$ and $1+2x\rho>0$ along the trajectory, so the sSQP subproblem has a unique global minimizer.

\noindent {\rm (ii)} Let $r_k:= \| (x_k, y_k, \lambda_k) \|.$ Equations \eqref{eq:t-recurrence} and \eqref{eq:iterates} give
\begin{equation*}
 r_{2j}=t_j(1+o(1)),\qquad
 r_{2j+1}=t_j^2(1+o(1)),\qquad
 r_{2j+2}=t_j^4(1+o(1)).
\end{equation*}
Therefore, $\frac{r_{2j+1}}{r_{2j}^2}\rightarrow1$ and $ \frac{r_{2j+2}}{r_{2j+1}^2}\rightarrow1$, and the proof of (ii) is complete.

\noindent {\rm (iii)} Since $y_k=0$, we have $e_k=|x_k|$. By~\eqref{eq:t-recurrence} and~\eqref{eq:iterates},
$$
 e_{2j}=t_j,\qquad
 e_{2j+1}=t_j^4(c+o(1)),\qquad
 e_{2j+2}=t_j^4(1+o(1)).
$$
Hence $e_{2j+1}/e_{2j}\to0$, whereas $e_{2j+2}/e_{2j+1}\to1/c=1+\sqrt{2}>1$.

\noindent {\rm (iv)} By the two limits in part~{\rm (iii)},
$$
 0\le \frac{e_{k+1}}{\max\{e_k,e_{k-1}\}}
 \le \frac{e_{k+1}}{e_{k-1}}
 =\frac{e_{k+1}}{e_k}\frac{e_k}{e_{k-1}}
 \rightarrow 0.
$$
\end{proof}

Along this trajectory, the primal error increases at every sufficiently large odd step, even though the primal-dual sequence converges $Q$-quadratically. We next examine which conditions on an SQP step guarantee primal $Q$-superlinear convergence.

\section{Characterizations for primal convergence}\label{sec:criteria}
Throughout this section, we assume that the primal-dual sequence $\{z_k\}$ converges to $\bar z$. We first study primal $Q$-superlinear convergence for sequences satisfying the pSQP equations~\eqref{eq:psqp} and then specialize the analysis to ordinary sSQP. We use the following notation.
\begin{align*}
 A_k& =A(x_k),\qquad G_k=G(x_k), \\ 
 s_k=\|d_k\|,\quad
 \tau_k&=\|(\Delta\lambda_k,\Delta\mu_k)\|, \quad e_k=\|x_k-\bar x\|,\quad
 \rho_k=\|\mathcal R(z_k)\|.
\end{align*}

\subsection{Primal superlinear convergence of perturbed SQP}

We begin with perturbations in the SQP equations and identify the components that determine the primal rate.  To measure stationarity in the critical directions, let $P_{\mathcal{C}}$ denote Euclidean projection onto the critical cone $\mathcal{C}$ defined in~\eqref{eq:critical-cone}. This gives the projected residual
\begin{equation}\label{eq:Theta}
 \Theta(x,\lambda,\mu) :=\|P_{\mathcal{C}}(\nabla_xL(x,\lambda,\mu))\| +\|h(x)\|+\|\min\{\mu,-g(x)\}\|.
\end{equation}
The following assumption records a projected error bound, which is satisfied under SOSC~\cite[Theorem~2.3]{fernandez2010sharp}.
\begin{assumption}[Projected error bound]\label{ass:projected-eb}
There exist a neighborhood $U$ of $\bar z$ and $\kappa>0$ such that
$$
 \|x-\bar x\|\leq\kappa\Theta(x,\lambda,\mu)
 \qquad\text{for every $z\in U$}.
$$
\end{assumption}

Consider a sequence $z_k\to\bar z$ generated by pSQP whose subproblem is specified in \eqref{eq:psqp}, and set
$$
 q_k:=-g(x_k)-G_kd_k.
$$
As in~\cite[Proposition~3.1]{fernandez2010sharp}, we assume that $(\omega_3^k)_{\mathcal N}\to0$. Under this assumption, the inactive multiplier components vanish for all sufficiently large $k$. The following lemma also gives an upper bound for $\Theta$ in terms of the primal error.

\begin{lemma}\label{lem:reverse}
Let $\mathcal{C}$ be the critical cone~\eqref{eq:critical-cone} at $\bar z$. The following assertions hold.
\begin{itemize}
\item[{\rm (i)}] There exist $\delta>0$ and $C_3>0$ such that $\Theta(x,\lambda,\mu)\le C_3\|x-\bar x\|$ for every $z=(x,\lambda,\mu)$ with $\|z-\bar z\|\le\delta$, $\mu\ge0$, and $\mu_{\mathcal N}=0$.

\item[{\rm (ii)}] If $z_k\to\bar z$ satisfies the pSQP equations~\eqref{eq:psqp} and $(\omega_3^k)_{\mathcal N}\to0$, then $\mu_{k+1,\mathcal N}=0$ for all sufficiently large $k$.
\end{itemize}
\end{lemma}
\begin{proof} \ 
\noindent {\rm (i)} Write $y=(\lambda,\mu)$ and fix $z=(x,y)$ satisfying the condition in statement. Since $\nabla_xL(\bar z)=0$ and $\mu_{\mathcal N}=0$, the definition of the critical cone in~\eqref{eq:critical-cone} implies that for every $v\in\mathcal C$,
$$
\begin{aligned}
 \langle\nabla_xL(\bar x,y),v\rangle
 & =\langle\nabla_xL(\bar x,y)-\nabla_xL(\bar z),v\rangle \\
 & = \langle \bar A^\top(\lambda-\bar\lambda)
  +\bar G_{\mathcal A_+}^\top(\mu_{\mathcal A_+}-\bar\mu_{\mathcal A_+})
  +\bar G_{\mathcal A_0}^\top\mu_{\mathcal A_0}, v\rangle  =\mu_{\mathcal A_0}^{\top}\bar G_{\mathcal A_0}v\le0,
 \end{aligned}
$$
where we use $\bar A v=0$, $\bar G_{\mathcal A_+}v=0$, and $\bar G_{\mathcal A_0}v\le0$ for $v\in \mathcal{C}$. By using projection theorem and $0\in\mathcal{C}$, we have $P_{\mathcal{C}}(\nabla_xL(\bar x,\lambda,\mu))=0$. Moreover, $h(\bar x)=0$, and $\mu\ge0$ and $\mu_{\mathcal N}=0$ imply $\min\{\mu,-g(\bar x)\}=0$. Thus $\Theta(\bar x,y)=0$.

By continuity, choose $\delta>0$ and $K>0$ such that $\|\nabla^2_{xx}L(x,y)\|$, $\|A(x)\|$, and $\|G(x)\|$ are bounded by $K$ whenever $\|z-\bar z\|\le\delta$. Using the nonexpansiveness of $P_{\mathcal C}$ and of the map $w\mapsto\min\{\mu,w\}$, we obtain
$$
\begin{aligned}
 \Theta(x,y)
 &\le \|\nabla_xL(x,y)-\nabla_xL(\bar x,y)\|
      +\|h(x)-h(\bar x)\|+\|g(x)-g(\bar x)\| \le 3K\|x-\bar x\|.
\end{aligned}
$$
This proves part~{\rm (i)} with $C_3=3K$.

\medskip\noindent {\rm (ii)} Since $z_k\to\bar z$, we have $d_k\to0$ and $G_k$ is bounded. Together with $(\omega_3^k)_{\mathcal N}\to0$, for each $i\in\mathcal N$,
$$
 q_{k,i}-(\omega_3^k)_i =-g_i(x_k)-g_i'(x_k)d_k-(\omega_3^k)_i \rightarrow-g_i(\bar x)>0.
$$
The complementarity equation in~\eqref{eq:psqp} therefore implies $\mu_{k+1,\mathcal N}=0$ for all sufficiently large $k$.
\end{proof}

Since $\mu_{k+1}\ge0$, Lemma~\ref{lem:reverse} and Assumption~\ref{ass:projected-eb} now imply, for all sufficiently large $k$,
$$
 \kappa^{-1}e_{k+1}\leq\Theta(z_{k+1})\leq C_3e_{k+1}.
$$
We can therefore study the primal error through $\Theta(z_{k+1})$. To relate this quantity to the pSQP equations, define
\begin{equation}\label{eq:chi-nu}
\begin{aligned}
 \chi_k&=\min\{\mu_{k+1},q_k\},\qquad \nu_k=\|P_{\mathcal{C}}(-\omega_1^k)\|+\|\omega_2^k\|+\|\chi_k\|.
\end{aligned}
\end{equation}
Since $\min\{\mu_{k+1,i},q_{k,i}-(\omega_3^k)_i\}=0$, the Lipschitz continuity of the scalar minimum gives
\begin{equation}\label{eq:chiomega}
 |\chi_{k,i}|\le |(\omega_3^k)_i|,\qquad i=1,\ldots,p.
\end{equation}
In particular, $\|\chi_k\|\le\|\omega_3^k\|$.

Note that $\nu_k$ and $\Theta(z_{k+1})$ measure projected stationarity, equality feasibility, and complementarity residuals at $z_{k+1}$ for the unperturbed SQP model and the original nonlinear problem, respectively. The next lemma shows that these two residual measures differ by $o(s_k)$.

\begin{lemma}\label{lem:transport}
Let $\{z_k\}$ satisfy the pSQP equations~\eqref{eq:psqp} and suppose that $z_k\to\bar z$. For $\Theta$ and $\nu_k$ defined in~\eqref{eq:Theta} and~\eqref{eq:chi-nu}, respectively, one has
\begin{equation}\label{eq:transport}
 |\Theta(z_{k+1})-\nu_k|=o(s_k).
\end{equation}
\end{lemma}
\begin{proof} \ 
Since $z_k\to\bar z$, both $s_k$ and $\tau_k$ tend to zero. From the first equation of~\eqref{eq:psqp},
\begin{equation}\label{eq:comp1bound}
\begin{aligned}
 \nabla_xL(z_{k+1})+\omega_1^k
 ={}&\nabla_xL(x_{k+1},\lambda_k,\mu_k)
       -\nabla_xL(z_k)-H_kd_k\\
    &+[A(x_{k+1})-A_k]^\top\Delta\lambda_k
      +[G(x_{k+1})-G_k]^\top\Delta\mu_k.
\end{aligned}
\end{equation}
The first line is $o(s_k)$ by continuity of the Hessian, and the second line is $O(s_k\tau_k)$, so
\begin{equation}\label{eq:comp1bound2}
 \|\nabla_xL(z_{k+1})+\omega_1^k\|=o(s_k).
\end{equation}
For the second pSQP equation, and the definition of $q_k$ give
\begin{equation}\label{eq:comp23bound}
\begin{aligned}
 h(x_{k+1})
 &=h(x_k)+A_kd_k+O(s_k^2)
 =-\omega_2^k+O(s_k^2),\\
 -g(x_{k+1})
 &=-g(x_k)-G_kd_k+O(s_k^2)
 =q_k+O(s_k^2).
\end{aligned}
\end{equation}
Combining the three components and using the nonexpansive property of the projection, we conclude that
$$
\begin{aligned}
 \bigl|\|P_{\mathcal{C}}(\nabla_xL(z_{k+1}))\|
          -\|P_{\mathcal{C}}(-\omega_1^k)\|\bigr|
 &\le\|\nabla_xL(z_{k+1})+\omega_1^k\|
 =o(s_k),\\
 \bigl|\|h(x_{k+1})\|-\|\omega_2^k\|\bigr|
 &\le\|h(x_{k+1})+\omega_2^k\|=O(s_k^2),\\
 \bigl|\|\min\{\mu_{k+1},-g(x_{k+1})\}\|-\|\chi_k\|\bigr|
 &\le\|-g(x_{k+1})-q_k\|=O(s_k^2).
\end{aligned}
$$
Since $s_k\to0$, the last two bounds are also $o(s_k)$. Adding the three bounds proves~\eqref{eq:transport}.
\end{proof}

We now use these estimates to characterize primal $Q$-superlinear convergence under either of the two error bounds.

\begin{theorem}\label{thm:effective}
Let $\{z_k\}$ be a sequence satisfying the pSQP equations~\eqref{eq:psqp} such that $z_k\to\bar z$. Assume $(\omega_3^k)_{\mathcal N}\to0$ and that either
\begin{itemize}
\item[{\rm (i)}] the projected error bound, i.e., Assumption~\ref{ass:projected-eb} holds; or
\item[{\rm (ii)}] the ordinary primal error bound, i.e., Assumption \ref{ass:eb} holds and $\|\omega_1^k\|=o(s_k)$.
\end{itemize}
Then
$$e_{k+1}=o(e_k) \quad \Longleftrightarrow \quad \nu_k=o(s_k).$$
\end{theorem}
\begin{proof} \ 
Since $|s_k-e_k|\le e_{k+1}$, we have
\begin{equation}\label{eq:primal-step-small-o}
 e_{k+1}=o(s_k)\quad\Longleftrightarrow\quad e_{k+1}=o(e_k).
\end{equation}
Indeed, if $e_{k+1}\le\varepsilon_k s_k$ for some $\varepsilon_k\ge0$ with $\varepsilon_k\to0$, then $(1-\varepsilon_k)e_{k+1}\le\varepsilon_k e_k.$ The converse follows by interchanging $s_k$ and $e_k$.

\medskip
\noindent\textbf{Sufficiency.} Suppose that $\nu_k=o(s_k)$. In case~{\rm (i)}, the projected error bound and Lemma~\ref{lem:transport} give
$$
 e_{k+1}\le\kappa\Theta(z_{k+1})
 \le\kappa\nu_k+o(s_k)=o(s_k).
$$
In case~{\rm (ii)}, \eqref{eq:comp1bound2} and $\|\omega_1^k\|=o(s_k)$ give $\|\nabla_xL(z_{k+1})\|=o(s_k)$. By~\eqref{eq:comp23bound}, the sum of the equality and complementarity residuals is at most $\|\omega_2^k\|+\|\chi_k\|+O(s_k^2)=o(s_k)$. Thus $\|\mathcal R(z_{k+1})\|=o(s_k)$, and Assumption~\ref{ass:eb} gives $e_{k+1}=o(s_k)$. In either case, \eqref{eq:primal-step-small-o} yields $e_{k+1}=o(e_k)$.

\medskip\noindent \textbf{Necessity.} Suppose that $e_{k+1}=o(e_k)$. By~\eqref{eq:primal-step-small-o}, $e_{k+1}=o(s_k)$. Lemma~\ref{lem:reverse}{\rm (ii)} gives $\mu_{k+1,\mathcal N}=0$ for all sufficiently large $k$, so Lemma~\ref{lem:reverse}{\rm (i)} yields $\Theta(z_{k+1})\le C_3e_{k+1}$. Hence, by Lemma~\ref{lem:transport}, $\nu_k\le\Theta(z_{k+1})+o(s_k)=o(s_k).$
\end{proof}

\begin{remark}\label{remark:psqp}
{\rm (i)} Suppose that SOSC holds at $\bar z$ and that a sequence satisfying~\eqref{eq:psqp} obeys $z_k\to\bar z$ and $(\omega_3^k)_{\mathcal N}\to0$. The necessary condition of primal superlinear convergence established in~\cite[Proposition~3.1 and Remark~3.1]{fernandez2010sharp} is
\begin{equation}\label{eq:fis-necessary}
 \|P_{\mathcal C}(-\omega_1^k)\|+\|\omega_2^k\|
 +\|(\omega_3^k)_{\mathcal A_+}\|=o(s_k).
\end{equation}
Under this condition, \cite[Theorem~3.4]{fernandez2010sharp} guarantees primal superlinear convergence if in addition $\|(\omega_3^k)_{\mathcal A_0}\|=o(s_k)$.  We take a closer look on weakly active indices. Under~\eqref{eq:fis-necessary} and \eqref{eq:chiomega}, we have $\|\chi_{k,\mathcal A_+}\|=o(s_k)$. For each $i\in\mathcal N$, Lemma~\ref{lem:reverse}{\rm (ii)} yields $\mu_{k+1,i}=0$ eventually, and $q_{k,i}\to-g_i(\bar x)>0$, so $\chi_{k,i}=0$ eventually. Thus $\nu_k=\|\chi_{k,\mathcal A_0}\|+o(s_k)$, and Theorem~\ref{thm:effective} yields
\begin{equation}\label{eq:fis-weak}
 e_{k+1}=o(e_k)
 \quad\Longleftrightarrow\quad
 \|\chi_{k,\mathcal A_0}\|=o(s_k).
\end{equation}

The condition $\|(\omega_3^k)_{\mathcal A_0}\|=o(s_k)$ is sufficient for~\eqref{eq:fis-weak}.  However, it may not be necessary. For example, if $\mu_{k+1,i}=0$ and $q_{k,i}\ge0$ at a weakly active index, then $\chi_{k,i}=0$, and the perturbed complementarity condition allows any $(\omega_3^k)_i\le0$.

\smallskip\noindent {\rm (ii)} Under the setting of Theorem~\ref{thm:effective}, suppose that either set of its assumptions holds and that, for all sufficiently large $k$ and each $i\in\mathcal A$,
\begin{equation}\label{eq:sign-compatible}
 \mu_{k+1,i}=0\quad\Longrightarrow\quad \omega_{3,i}^k\geq0.
\end{equation}
Then $\|(\omega_3^k)_{\mathcal A}\|=o(s_k)$ is also necessary for primal superlinear convergence. More precisely,
\begin{equation}\label{eq:sign-iff}
 e_{k+1}=o(e_k)\quad\Longleftrightarrow\quad
 \|P_{\mathcal{C}}(-\omega_1^k)\|+\|\omega_2^k\|
     +\|(\omega_3^k)_{\mathcal A}\|=o(s_k).
\end{equation}
Indeed, if $e_{k+1}=o(e_k)$, then $e_{k+1}=o(s_k)$ by~\eqref{eq:primal-step-small-o}. Since $g_{\mathcal A}(\bar x)=0$, \eqref{eq:comp23bound} gives
$$
 \|q_{k,\mathcal A}\|  \le \|g_{\mathcal A}(x_{k+1})\|+O(s_k^2)  = O(e_{k+1})+O(s_k^2)=o(s_k).
$$
At an active index with $\mu_{k+1,i}>0$, complementarity gives $\omega_{3,i}^k=q_{k,i}$. If $\mu_{k+1,i}=0$, feasibility and~\eqref{eq:sign-compatible} give $0\le\omega_{3,i}^k\le q_{k,i}$. Thus $(\omega_3^k)_{\mathcal A}=o(s_k)$, and the other two terms can be derived by Theorem~\ref{thm:effective}. Conversely, \eqref{eq:chiomega} gives $\|\chi_{k,\mathcal A}\|\le\|(\omega_3^k)_{\mathcal A}\|$. Since $\chi_{k,\mathcal N}=0$ eventually, the right-hand side of~\eqref{eq:sign-iff} implies $\nu_k=o(s_k)$. Theorem~\ref{thm:effective} yields primal $Q$-superlinear convergence.
\end{remark}

\subsection{Primal superlinear convergence of sSQP}
We now apply Theorem~\ref{thm:effective} to the sSQP iteration~\eqref{eq:ssqp-system}, which satisfies~\eqref{eq:psqp} with the perturbations~\eqref{eq:ssqp-omega}. In particular, the stationarity perturbation vanishes, so case {\rm (ii)} of Theorem~\ref{thm:effective} applies under Assumption \ref{ass:eb}. The following theorem gives equivalent criteria for primal $Q$-superlinear convergence in terms of the primal and multiplier steps and the KKT residual.

\begin{theorem}\label{thm:product}
Assume Assumption \ref{ass:eb} holds. Let $\{z_k\}_{k\in\mathbb{N}}$ be a sequence generated by~\eqref{eq:ssqp-system} and satisfy  $z_k\to\bar z$. Then the following four statements are equivalent:
\begin{equation*}
 e_{k+1}=o(e_k),\qquad
 \rho_k\tau_k=o(s_k),\qquad
 \tau_k^2=o(s_k), \qquad
 \rho_k^2=o(e_k).
\end{equation*}
\end{theorem}

\begin{proof}\ 
By~\eqref{eq:ssqp-omega}, $\omega_1^k=0$ and $\|(\omega_3^k)_{\mathcal N}\|\le\rho_k\tau_k\to0$. If $\mu_{k+1,i}=0$, then $\omega_{3,i}^k=\rho_k\mu_{k,i}\ge0$, so~\eqref{eq:sign-compatible} holds. Moreover, Lemma~\ref{lem:reverse}{\rm (ii)} gives $\Delta\mu_{k,\mathcal N}=0$ for all sufficiently large $k$. Consequently, it follows from Remark~\ref{remark:psqp}{\rm (ii)} and \eqref{eq:ssqp-omega} that
$$
 e_{k+1}=o(e_k)\quad\Longleftrightarrow\quad  \rho_k\bigl(\|\Delta\lambda_k\|  +\|\Delta\mu_{k,\mathcal A}\|\bigr)=o(s_k) \quad\Longleftrightarrow\quad
 \rho_k\tau_k=o(s_k),
$$
where the second equivalence follows from $\tau_k \le\|\Delta\lambda_k\|+\|\Delta\mu_{k,\mathcal A}\| \le\sqrt2\,\tau_k.$

By Proposition~\ref{prop:local-solvability}{\rm (i)} and Lemma~\ref{lem:upper}, we have $\tau_k\le C_1\rho_k$ and $\rho_k\le C_2(s_k+\tau_k)$. Hence
$$
 \tau_k^2\le C_1\rho_k\tau_k,\qquad \rho_k\tau_k\le C_2(s_k\tau_k+\tau_k^2).
$$
Since $\tau_k\to0$, these inequalities show that
$$\rho_k\tau_k=o(s_k) \quad\Longleftrightarrow\quad \tau_k^2=o(s_k).$$

Finally, suppose that $e_{k+1}=o(e_k)$. We have from the triangle inequality that $s_k\le e_k+e_{k+1}=O(e_k)$. From the equivalence established above, we have $\tau_k^2=o(s_k)$. Hence $\rho_k^2\le 2C_2^2(s_k^2+\tau_k^2)=o(e_k)$. Conversely, suppose that $\rho_k^2=o(e_k)$. By~\eqref{eq:comp1bound2}--\eqref{eq:comp23bound}, \eqref{eq:ssqp-omega}, and~\eqref{eq:chiomega},
$$
 \rho_{k+1}\le o(s_k)+O(s_k^2+\rho_k\tau_k).
$$
Using Assumption~\ref{ass:eb} and $\tau_k\le C_1\rho_k$, we obtain $e_{k+1}\le\varepsilon_k s_k+O(\rho_k^2)$ for some $\varepsilon_k\ge0$ with $\varepsilon_k\to0$. Since $s_k\le e_k+e_{k+1}$, it follows that
$$
 (1-\varepsilon_k)e_{k+1}
 \le\varepsilon_k e_k+O(\rho_k^2)=o(e_k).
$$
This proves the remaining equivalence.
\end{proof}

Next, we revisit the counterexample in Section~\ref{sec:counterexample} to show that the condition $\tau_k^2=o(s_k)$ in Theorem~\ref{thm:product} cannot in general be weakened to $\tau_k^2=O(s_k)$.

\begin{corollary}\label{cor:threshold}
Along the sequence specified in Theorem~\ref{thm:counterexample},
$$
 \frac{\tau_{2j+1}^2}{s_{2j+1}}\rightarrow\frac1{\sqrt2},
 \qquad
 \frac{\rho_{2j+1}\tau_{2j+1}}{s_{2j+1}}
 \rightarrow\frac1{\sqrt2}.
$$
Consequently, $\tau_k^2=O(s_k)$ along the full trajectory, and $\tau_k^2=o(s_k)$ fails.
\end{corollary}
\begin{proof} \ 
Equations~\eqref{eq:t-recurrence} and~\eqref{eq:iterates}, together with the definition of $\rho_k$, give
$$
 \frac{s_{2j+1}}{t_j^4}\to\sqrt2,\qquad
 \frac{\tau_{2j+1}}{t_j^2}\to1,\qquad
 \frac{\rho_{2j+1}}{t_j^2}\to1.
$$
Consequently,
$$
 \frac{\tau_{2j+1}^2}{s_{2j+1}}\to\frac1{\sqrt2},
 \qquad
 \frac{\rho_{2j+1}\tau_{2j+1}}{s_{2j+1}}\to\frac1{\sqrt2}.
$$
Thus, we have that $\tau_k^2=o(s_k)$ fails.

For even indices, \eqref{eq:iterates} gives $s_{2j}=t_j(1+o(1))$ and $\tau_{2j}=O(t_j^3)$, so
$$
 \frac{\tau_{2j}^2}{s_{2j}}=O(t_j^5)\to0.
$$
Thus $\tau_k^2/s_k$ is bounded along both the even and odd indices, proving $\tau_k^2=O(s_k)$.
\end{proof}

We next relate the current KKT residual to the primal errors before and after an sSQP step.

\begin{lemma}\label{lem:local-primal-scale}
Suppose that Assumptions~\ref{ass:eb} and~\ref{ass:lipschitz} hold. Write $e:=\|x-\bar x\|$ and $e^+:=\|x^+-\bar x\|$. Then there exist a neighborhood $U$ of $\bar z$ and a constant $C_4\ge1$ such that the following assertions hold for every pair $z,z^+\in U$ satisfying~\eqref{eq:ssqp-system} with $\mu\ge0$ and $\rho>0$.

\begin{itemize}
\item[{\rm (i)}] We have $\mu^+_{\mathcal N}=0$ and $ \|\mathcal R(z^+)\|+e^+\le C_4\rho^2.$

\item[{\rm (ii)}] If $\mu_{\mathcal N}=0$, then $C_4^{-1}\rho^2\le e^2+e^+\le C_4\rho^2.$
\end{itemize}
\end{lemma}

\begin{proof} \ 
Let $d:=x^+-x$, $s:=\|d\|$, and $\tau:=\|(\Delta\lambda,\Delta\mu)\|$. By Proposition~\ref{prop:local-solvability}{\rm (i)} and Lemma~\ref{lem:upper},
\begin{equation}\label{eq:rhostau}
 C_2^{-1}\rho\le (s^2+\tau^2)^{1/2}\le C_1\rho.
\end{equation}

\medskip\noindent {\rm (i)} For $i\in\mathcal N$, the linearized slack in~\eqref{eq:ssqp-system} tends to $-g_i(\bar x)>0$ as $z,z^+\to\bar z$. Hence $\mu^+_{\mathcal N}=0$ by complementarity. Under Assumption~\ref{ass:lipschitz}, it follows from the expression in \eqref{eq:comp1bound} that $\|\nabla_xL(z^+)\|=O(s^2+s\tau)$. Using \eqref{eq:comp23bound}, \eqref{eq:ssqp-omega}, and the nonexpansiveness of the componentwise minimum, we obtain
$$
 \|h(x^+)\|+\|\min\{\mu^+,-g(x^+)\}\|=O(s^2+\rho\tau).
$$
Since $\max\{s,\tau\}\le C_1\rho$, we obtain $\|\mathcal R(z^+)\|=O(\rho^2)$. Then Assumption~\ref{ass:eb} implies $e^+\le\kappa\|\mathcal R(z^+)\|=O(\rho^2)$.

\medskip\noindent {\rm (ii)} By part~{\rm (i)} and $e\le\kappa\rho$, we have $e^2+e^+=O(\rho^2)$. To prove the lower bound, set $F=(h,g_{\mathcal A})$. Since $\mu_{\mathcal N}=0$ by assumption and $\mu^+_{\mathcal N}=0$ by part~{\rm (i)}, we have $\Delta\mu_{\mathcal N}=0$. The sSQP equations are formulated as
$$
 \rho\Delta\lambda=h(x)+A(x)d,\qquad  \mu^+=\bigl[\mu+(g(x)+G(x)d)/\rho\bigr]_+.
$$
Since $\mu\ge0$, the Lipschitz continuity of $[\cdot]_+$ implies $\rho\|\Delta\mu_{\mathcal A}\|\le\|g_{\mathcal A}(x)+G_{\mathcal A}(x)d\|$. Together with the equality above, this yields $\rho\tau\le\|F(x)+F'(x)d\|$. Expanding $F(\bar x)=0$ at $x$, we have
$$
 F(x)+F'(x)d=F'(x)(x^+-\bar x)+O(e^2).
$$
Since $\tau\le C_1\rho$ by Proposition~\ref{prop:local-solvability}{\rm (i)}, we obtain
$$
 \tau^2\le C_1\rho\tau  \le C_1\|F(x)+F'(x)d\| =O(e^2+e^+).
$$
Using $s\le e+e^+$ and taking $U$ small enough that $e^+\le1$, we obtain from \eqref{eq:rhostau} that
$$
 \rho^2\le C_2^2(s^2+\tau^2)=O(e^2+e^+).
$$
The estimates are uniform for $z,z^+\in U$. Taking $C_4$ sufficiently large gives both assertions.
\end{proof}

We now refine the criterion for primal $Q$-superlinear convergence by characterizing the rates $e_{k+1}=O(e_k^\gamma)$ for $1<\gamma\le2$ in terms of the KKT residual and the step lengths.

\begin{theorem}[Primal convergence rates]\label{thm:primal-order}
Suppose that the assumptions of Theorem~\ref{thm:product} and Assumption~\ref{ass:lipschitz} hold. For every $1<\gamma\le2$,
$$
 e_{k+1}=O(e_k^\gamma)
 \quad\Longleftrightarrow\quad \rho_k^2=O(e_k^\gamma)
 \quad\Longleftrightarrow\quad \tau_k^2=O(s_k^\gamma).
$$
In particular, the following conditions are equivalent:
\begin{equation}\label{eq:primal-quadratic-criterion}
 e_{k+1}=O(e_k^2),\qquad \rho_k=O(e_k),\qquad
 \rho_k=O(s_k),\qquad \tau_k=O(s_k).
\end{equation}
\end{theorem}

\begin{proof} \ 
Since $e_k^2=O(e_k^\gamma)$, the first equivalence follows from Lemma~\ref{lem:local-primal-scale}{\rm (ii)}. By Proposition~\ref{prop:local-solvability}{\rm (i)}, $\tau_k\le C_1\rho_k$. Also, Lemma~\ref{lem:local-primal-scale}{\rm (i)} and Lemma~\ref{lem:upper} give
\begin{equation}\label{eq:quantitative-bound}
 e_{k+1}\le C_4\rho_k^2\le 2C_4C_2^2(s_k^2+\tau_k^2).
\end{equation}

If $e_{k+1}=O(e_k^\gamma)$, then $e_{k+1}\le e_k/2$. Since $|s_k-e_k|\le e_{k+1}$, we have $e_k\le2s_k$, and therefore $\tau_k^2=O(\rho_k^2)=O(e_k^\gamma)=O(s_k^\gamma)$. Conversely, if $\tau_k^2=O(s_k^\gamma)$, then \eqref{eq:quantitative-bound} gives $e_{k+1}=O(s_k^\gamma)$. As $\gamma>1$, we have $e_{k+1}\le s_k/2$ eventually. It follows that $s_k\le2e_k$ and hence $e_{k+1}=O(e_k^\gamma)$. Taking $\gamma=2$ and using the same step bounds yields \eqref{eq:primal-quadratic-criterion}.
\end{proof}

\section{Multiplier correction for stabilized SQP}\label{relax:section}

In the counterexample of Section~\ref{sec:counterexample}, the ratio $\rho_k/e_k$ tends to infinity along the odd iterates, and the primal error increases at the following step. For a convergent ordinary sSQP sequence, Theorem~\ref{thm:primal-order} characterizes primal $Q$-quadratic convergence by $\rho_k=O(e_k)$ and shows that this condition is equivalent to $\rho_k=O(s_k)$. Both $\rho_k$ and $s_k$ are available once the sSQP subproblem has been solved.

This characterization motivates the multiplier correction in Algorithm~\ref{alg:corrected-ssqp}. We retain the primal point $\widehat x_k$ computed by ordinary sSQP and compare $\rho_k$ with $b_k^2s_k$ to choose the weight on the current multipliers in the update. Our aim is to obtain $\rho_k=O(e_k)$ eventually along the corrected sequence. Indeed, since $x_{k+1}=\widehat x_k$, Lemma~\ref{lem:local-primal-scale}{\rm (i)} gives $e_{k+1}=O(\rho_k^2)$, so the desired bound on $\rho_k$ yields $e_{k+1}=O(e_k^2)$.

Throughout this section, we assume that Assumption~\ref{ass:lipschitz} holds and use the multiplier set $\mathcal M_*$ defined in~\eqref{eq:defM}.

\begin{algorithm}[!ht]
\caption{Stabilized SQP with multiplier correction}
\label{alg:corrected-ssqp}
\textbf{Initialization.} Fix a sufficiently large constant $\ell_0>0$. Choose $b_0>0$ and $z_0:=(x_0,\lambda_0,\mu_0)$ with $\mu_0\ge0$. Set $k:=0$.

\smallskip\noindent \textbf{Step 1:} Let $\rho_k:=\|\mathcal R(z_k)\|$. If $\rho_k=0$, stop; otherwise compute $\widehat z_k:=(\widehat x_k,\widehat\lambda_k,\widehat\mu_k)$ satisfying
\begin{equation}\label{relax:predictor}
\begin{aligned}
0={}&\nabla_xL(z_k)+H_kd_k
       +A_k^\top\Delta\widehat\lambda_k
       +G_k^\top\Delta\widehat\mu_k,\\
0={}&h(x_k)+A_kd_k-\rho_k\Delta\widehat\lambda_k,\\
0\le{}&\widehat\mu_k\ \perp\
       -g(x_k)-G_kd_k+\rho_k\Delta\widehat\mu_k\ge0,
\end{aligned}
\end{equation}
and
\begin{equation}\label{relax:localization}
 \|(d_k,\Delta\widehat y_k)\|\le\ell_0\rho_k,
\end{equation}
where
$$
\begin{gathered}
 H_k:=\nabla^2_{xx}L(z_k),\quad d_k:=\widehat x_k-x_k,\\
 y_k:=(\lambda_k,\mu_k),\quad
 \widehat y_k:=(\widehat\lambda_k,\widehat\mu_k), \quad
 \Delta\widehat y_k
 :=(\Delta\widehat\lambda_k,\Delta\widehat\mu_k)
 =\widehat y_k-y_k.
\end{gathered}
$$
\smallskip\noindent\textbf{Step 2:} Let $s_k := \|d_k\|$ and set
\begin{equation}\label{relax:adaptive-test}
 c_k:=\begin{cases}
  2b_k,&\rho_k>b_k^2s_k,\\
  b_k,&\rho_k\le b_k^2s_k,
 \end{cases}
 \quad b_{k+1}:=c_k,\quad
 \omega_k:=\min\{1/2,c_k\rho_k\}.
\end{equation}
Define the index set
\begin{equation*}
 \mathcal I_k:=\{i:g_i(x_k)\ge-\sqrt{\rho_k}\},
\end{equation*}
and set
\begin{equation}\label{relax:update}
\begin{aligned}
 x_{k+1}&:=\widehat x_k,\quad
 \lambda_{k+1}:=\omega_k\lambda_k+(1-\omega_k)\widehat\lambda_k,\\
 \mu_{k+1,i}&:=\begin{cases}
 \omega_k\mu_{k,i}+(1-\omega_k)\widehat\mu_{k,i},&i\in\mathcal I_k,\\
 0,&i\notin\mathcal I_k.
 \end{cases}
\end{aligned}
\end{equation}

Set $k\leftarrow k+1$ and return to Step~1.
\end{algorithm}

\begin{remark}
{\rm (i)} Under SOSC, Proposition~\ref{prop:local-solvability}{\rm (ii)} ensures the existence of a solution to~\eqref{relax:predictor} satisfying~\eqref{relax:localization} for every sufficiently close input with $\mu_k\ge0$ and $\rho_k>0$, provided that $\ell_0$ is sufficiently large. The condition~\eqref{relax:localization} restricts the choice of the subproblem solution to a neighborhood of the current iterate. Similar requirements appear in local convergence analysis; see~\cite[Theorem~1]{fernandez2014stabilized} and~\cite[Theorem~2]{izmailov2015combining}.

\smallskip\noindent {\rm (ii)} The set $\mathcal I_k$ estimates which inequalities are active at the solution. For $i\notin\mathcal I_k$, the update sets $\mu_{k+1,i}=0$, consistent with the complementarity condition at the solution.  The bound $\omega_k\le1/2$ in~\eqref{relax:adaptive-test} ensures that the weight on the predicted multipliers remains at least $1/2$, even when $c_k\rho_k$ is large. In particular, for the equality multipliers,
$$
\lambda_{k+1}-\lambda_k =(1-\omega_k)\Delta\widehat\lambda_k,\qquad \tfrac12\le1-\omega_k\le1,
$$
and this also holds for the inequality multipliers with $i\in\mathcal I_k$.

\smallskip \noindent {\rm (iii)} Multiplier adjustment has also been used to improve local convergence. Fischer's method~\cite{fischer1999modified} computes multiplier estimates from an auxiliary quadratic program, while Wright~\cite{wright2003constraint} uses an auxiliary linear program to adjust the estimates on selected constraints. Our method retains the ordinary sSQP primal point and adjusts the multipliers using only elementary vector and scalar operations, without solving an additional optimization subproblem.
\end{remark}

To analyze the multiplier correction, we use $\nabla_xL(\bar x,y_k)$ to measure the effect of the current multiplier error on stationarity, with the primal variable fixed at its solution $\bar x$. The following lemma establishes local estimates for the ordinary sSQP step and the KKT residual.

\begin{lemma}\label{relax:mixed-predictor}
Suppose that Assumptions~\ref{ass:eb} and~\ref{ass:lipschitz} hold. There exist $\delta>0$ and $C_5>0$ such that, for every input $z_k$ with $\|z_k-\bar z\|\le\delta$ and $\mu_k\ge0$, the following assertions hold.
\begin{itemize}
\item[{\rm (i)}] If $\rho_k>0$, every $\widehat z_k$ satisfying~\eqref{relax:predictor}--\eqref{relax:localization} obeys $s_k+\|\Delta\widehat y_k\|\le C_5\rho_k$, $\widehat\mu_{k,\mathcal N}=0$, and $\mathcal I_k=\mathcal A$.
\item[{\rm (ii)}] Every $\widehat z_k$ in {\rm (i)} also satisfies $\|\mathcal R(\widehat z_k)\|\le C_5\rho_k^2$ and $\|\widehat x_k-\bar x\|+\|\nabla_xL(\bar x,\widehat y_k)\|\le C_5\rho_k^2$.
\item[{\rm (iii)}] If $\mu_{k,\mathcal N}=0$, then $C_5^{-1}\rho_k\le\|x_k-\bar x\|+\|\nabla_xL(\bar x,y_k)\|\le C_5\rho_k.$
\end{itemize}
\end{lemma}
\begin{proof}\ 
Choose $\bar\delta> 0$ such that Assumptions~\ref{ass:eb} and \ref{ass:lipschitz} hold for $\|z-\bar z\|\le\bar\delta$, and let $M$ be the constant in~\eqref{eq:local-smoothness-bounds}. Choose $0<\delta\le\bar\delta/2$ with $\ell_0M\delta\le\bar\delta/2$. For every input satisfying $\|z_k-\bar z\|\le\delta$, it follows from \eqref{eq:local-smoothness-bounds} and $\mathcal R(\bar z)=0$ that $\rho_k\le M\|z_k-\bar z\|$. Hence, by~\eqref{relax:localization},
\begin{equation}\label{eq:zkhat}
 \|\widehat z_k-\bar z\| \le \|z_k-\bar z\|+\ell_0\rho_k \le (1+\ell_0M)\delta \le\bar\delta.
\end{equation}
For parts~{\rm (i)} and~{\rm (ii)}, suppose that $\rho_k>0$. By~\eqref{eq:zkhat}, after reducing $\delta$ if necessary, we have from Lemma~\ref{lem:local-primal-scale}{\rm (i)} that
\begin{equation}\label{eq:Rkxkhatbound}
 \widehat\mu_{k,\mathcal N}=0,\qquad
 \|\mathcal R(\widehat z_k)\|+\|\widehat x_k-\bar x\|
 \le C_4\rho_k^2.
\end{equation}

\medskip\noindent {\rm (i)} By~\eqref{relax:localization}, $s_k+\|\Delta\widehat y_k\|\le 2\ell_0\rho_k$. For $i\in\mathcal A$, we have $g_i(\bar x)=0$, which together with \eqref{eq:local-smoothness-bounds} and Assumptions~\ref{ass:eb} yields
$$
 |g_i(x_k)|\le M e_k\le M\kappa\rho_k\le\sqrt{\rho_k}
$$
for $\delta$ sufficiently small. Hence $i\in\mathcal I_k$ and $\mathcal A\subseteq\mathcal I_k$. For $i\in\mathcal N$, we have $g_i(x_k)\to g_i(\bar x)<0$ as $z_k\to\bar z$, while $\rho_k\to0$. Thus $g_i(x_k)<-\sqrt{\rho_k}$ for $z_k$ sufficiently close to $\bar z$. It follows that $\mathcal N\cap\mathcal I_k=\varnothing$, so $\mathcal I_k\subseteq\mathcal A$. Combining these two inclusions, we have $\mathcal I_k=\mathcal A$.

\medskip\noindent {\rm (ii)} By~\eqref{eq:local-smoothness-bounds}, we have $\|\nabla_xL(\bar x,\widehat y_k)\|  \le \|\mathcal R(\widehat z_k)\|+M\|\widehat x_k-\bar x\|.$ Adding $\|\widehat x_k-\bar x\|$ and using equation \eqref{eq:Rkxkhatbound} yield
$$
 \|\widehat x_k-\bar x\|+\|\nabla_xL(\bar x,\widehat y_k)\|
 \le (1+M)C_4\rho_k^2.
$$

\medskip\noindent {\rm (iii)} Since $h(\bar x)=0$, $\mu_k\ge0$, and $\mu_{k,\mathcal N}=0$, we have $\mathcal R(\bar x,y_k)=(\nabla_xL(\bar x,y_k),0,0)$. By the reverse triangle inequality and~\eqref{eq:local-smoothness-bounds},
$$
 \bigl|\rho_k-\|\nabla_xL(\bar x,y_k)\|\bigr|
 \le \|\mathcal R(x_k,y_k)-\mathcal R(\bar x,y_k)\|
 \le M e_k.
$$
Using $e_k\le\kappa\rho_k$ and $M\ge1$, we obtain
$$
 M^{-1}\rho_k
 \le e_k+\|\nabla_xL(\bar x,y_k)\|
 \le \bigl[1+\kappa(1+M)\bigr]\rho_k.
$$
Thus all three assertions hold with $C_5=\max\bigl\{2\ell_0,\,(1+M)C_4,\,M,\,1+\kappa(1+M)\bigr\}.$
\end{proof}

Since $x_{k+1}=\widehat x_k$, Lemma~\ref{relax:mixed-predictor}{\rm (ii)} implies $e_{k+1}\le C_5\rho_k^2$. This does not yet imply a quadratic primal rate, because $\rho_k/e_k$ may be large. Part~{\rm (iii)} shows that the difficult case occurs when $e_k$ is small relative to $\|\nabla_xL(\bar x,y_k)\|$. The next lemma considers this case and shows that $e_{k+1}$ is bounded below by a positive multiple of $\rho_k^2$.

\begin{lemma}\label{relax:mixed-lower}
Suppose that Assumptions~\ref{ass:eb} and~\ref{ass:lipschitz} hold. There exist $\varepsilon>0$, $\delta>0$, and $C_6>0$ such that, for every input $z_k$ with $\|z_k-\bar z\|\le\delta$, $\rho_k>0$, $\mu_k\ge0$, $\mu_{k,\mathcal N}=0$, and
\begin{equation}\label{eq:rel}
 \|x_k-\bar x\|\le\varepsilon\|\nabla_xL(\bar x,y_k)\|,
\end{equation}
every point $\widehat{z}_k$ satisfying~\eqref{relax:predictor}-\eqref{relax:localization} obeys
\begin{equation}\label{relax:mixed-cone-estimates}
 C_6^{-1}\rho_k^2\le\|\widehat x_k-\bar x\|\le C_6\rho_k^2.
\end{equation}
\end{lemma}
\begin{proof} \ 
Set $\widehat e_k:=\|\widehat x_k-\bar x\|$. By~\eqref{relax:localization}, both $z_k$ and $\widehat z_k$ are sufficiently close to $\bar z$ for small $\delta$. Since $\mu_{k,\mathcal N}=0$, from Lemma~\ref{lem:local-primal-scale}{\rm (ii)} we have
\begin{equation}\label{eq:rhokekehat}
 C_4^{-1}\rho_k^2\le e_k^2+\widehat e_k\le C_4\rho_k^2.
\end{equation}
Choose $0<\varepsilon\le1$ so that $\varepsilon^2C_5^2\le 1/(2C_4)$. It follows from \eqref{eq:rel} and Lemma~\ref{relax:mixed-predictor}{\rm (iii)} that
$$
 e_k\le\varepsilon\|\nabla_xL(\bar x,y_k)\| \le\varepsilon C_5\rho_k.
$$
Hence $e_k^2\le\rho_k^2/(2C_4)$. Substituting this bound into \eqref{eq:rhokekehat}, we obtain $\widehat e_k\ge C_4^{-1}\rho_k^2-e_k^2  \ge \frac{\rho_k^2}{2C_4}.$ The upper bound $\widehat e_k\le C_4\rho_k^2$ follows directly from $e_k^2+\widehat e_k\le C_4\rho_k^2$. Together, these bounds give~\eqref{relax:mixed-cone-estimates} with $C_6=2C_4$.
\end{proof}

We now explain the role of retaining part of the current multipliers. For $\mu_{k,\mathcal N}=0$, Lemma~\ref{relax:mixed-predictor}{\rm (i)} and \eqref{relax:update} imply that when $\omega_k=c_k\rho_k<1/2$,
$$
 \nabla_xL(\bar x,y_{k+1}) =c_k\rho_k\nabla_xL(\bar x,y_k) +(1-c_k\rho_k)\nabla_xL(\bar x,\widehat y_k).
$$
Under~\eqref{eq:rel}, Lemma~\ref{relax:mixed-predictor}{\rm (iii)} shows that the first term has magnitude proportional to $c_k\rho_k^2$, while part~{\rm (ii)} bounds the second term by $C_5\rho_k^2$. When $\rho_k$ is small, $c_k$ can be chosen large enough that the second term cannot nearly cancel the first, while still satisfying $c_k\rho_k<1/2$. Since $e_{k+1}\le C_5\rho_k^2$ as well, both the second term and the primal error are small compared with the first term when $c_k$ is sufficiently large. This suggests that retaining part of the current multipliers may help preserve condition~\eqref{eq:rel} at the next iterate. The next lemma proves this claim and gives an upper bound for $\rho_{k+1}$.

\begin{lemma}\label{relax:corrected-step}
Suppose that Assumptions~\ref{ass:eb} and~\ref{ass:lipschitz} hold. There exist $\delta>0$ and $C_7>0$ such that the following estimates hold for every input $z_k$ with $\|z_k-\bar z\|\le\delta$, $\rho_k>0$, $\mu_k\ge0$, and $\mu_{k,\mathcal N}=0$, every vector $\widehat{z}_k$ satisfying~\eqref{relax:predictor}--\eqref{relax:localization}, and the correction~\eqref{relax:update} with $0<\omega_k=c_k\rho_k<1/2$:
\begin{equation}\label{relax:accepted-estimates}
\begin{aligned}
 \|x_{k+1}-\bar x\|&\le C_7\rho_k^2,\\
 \|\nabla_xL(\bar x,y_{k+1})\|
 &\le c_k\rho_k\|\nabla_xL(\bar x,y_k)\|
       +C_7\rho_k^2,\\
 \rho_{k+1}&\le C_7(1+c_k)\rho_k^2.
\end{aligned}
\end{equation}
Moreover, for each fixed $\varepsilon\in(0,1]$, there exists a threshold $c_*>1$ such that, if
$$
 \|x_k-\bar x\|\le\varepsilon\|\nabla_xL(\bar x,y_k)\|
 \quad {\rm and}\quad c_k\ge c_*,
$$
then
\begin{equation}\label{relax:cone-invariance}
 \|x_{k+1}-\bar x\|
 \le\frac{\varepsilon}{2}\|\nabla_xL(\bar x,y_{k+1})\|.
\end{equation}
\end{lemma}
\begin{proof} \ 
Fix $C\ge\max\{1,C_5\}$, where $C_5$ is the constant in Lemma~\ref{relax:mixed-predictor}. By Lemma~\ref{relax:mixed-predictor}{\rm (i)}, we have $\mathcal I_k=\mathcal A$ and $\widehat\mu_{k,\mathcal N}=0$. Since also $\mu_{k,\mathcal N}=0$, it follows from~\eqref{relax:update} that, for $0\le\omega_k\le1/2$,
$$
 \mu_{k+1,\mathcal N} =\omega_k\mu_{k,\mathcal N} +(1-\omega_k)\widehat\mu_{k,\mathcal N}=0.
$$
Thus~\eqref{relax:update} can be written as
\begin{equation}\label{eq:yk+1ykhat}
 y_{k+1}=\omega_k y_k+(1-\omega_k)\widehat y_k.
\end{equation}
In particular, $\mu_{k+1}\ge0$ and $\mu_{k+1,\mathcal N}=0$. By~\eqref{relax:update}, \eqref{relax:localization}, and~\eqref{eq:local-smoothness-bounds},
$$
\begin{aligned}
 \|z_{k+1}-\bar z\|
 &\le \|z_k-\bar z\|
      +\|(d_k,(1-\omega_k)\Delta\widehat y_k)\| \le \|z_k-\bar z\|+\ell_0\rho_k \le (1+\ell_0M)\|z_k-\bar z\|.
\end{aligned}
$$
By reducing $\delta$ if necessary, we may apply Lemma~\ref{relax:mixed-predictor}{\rm (iii)} at $z_{k+1}$, since $\mu_{k+1}\ge0$ and $\mu_{k+1,\mathcal N}=0$.

Since the Lagrangian gradient is affine in the multipliers, we also have
\begin{equation}\label{relax:eta-recurrence}
\nabla_xL(\bar x,y_{k+1})
 =\omega_k\nabla_xL(\bar x,y_k)
  +(1-\omega_k)\nabla_xL(\bar x,\widehat y_k).
\end{equation}
For $\omega_k=c_k\rho_k$, it follows from Lemma~\ref{relax:mixed-predictor}{\rm (ii)} and~\eqref{relax:eta-recurrence} that
$$
\begin{aligned}
 \|x_{k+1}-\bar x\|
 &=\|\widehat x_k-\bar x\|\le C_5\rho_k^2,\\
 \|\nabla_xL(\bar x,y_{k+1})\|
 &\le c_k\rho_k\|\nabla_xL(\bar x,y_k)\|
       +C_5\rho_k^2.
\end{aligned}
$$
Applying the pointwise residual comparison in Lemma~\ref{relax:mixed-predictor}{\rm (iii)} at $z_k$ and $z_{k+1}$, we have
$$
 \|\nabla_xL(\bar x,y_k)\|\le C_5\rho_k,
 \qquad
 \rho_{k+1}\le C_5\bigl(\|x_{k+1}-\bar x\|
       +\|\nabla_xL(\bar x,y_{k+1})\|\bigr).
$$
Since $C_5\le C$, these estimates imply 
$$
 \rho_{k+1}
 \le2C^2\rho_k^2+C c_k\rho_k\|\nabla_xL(\bar x,y_k)\|
 \le C^2(2+c_k)\rho_k^2.
$$
Thus~\eqref{relax:accepted-estimates} holds with $C_7=2C^2$.

Now fix $\varepsilon\in(0,1]$ and suppose that $\|x_k-\bar x\|\le\varepsilon\|\nabla_xL(\bar x,y_k)\|$. By Lemma~\ref{relax:mixed-predictor}{\rm (iii)}, $\rho_k\le2C_5\|\nabla_xL(\bar x,y_k)\|\le2C\|\nabla_xL(\bar x,y_k)\|$. Using this bound, Lemma~\ref{relax:mixed-predictor}{\rm (ii)} and \eqref{relax:eta-recurrence}, we obtain
$$
\begin{aligned}
 \|x_{k+1}-\bar x\|
 &\le C\rho_k^2
 \le2C^2\rho_k\|\nabla_xL(\bar x,y_k)\|,\\
 \|\nabla_xL(\bar x,y_{k+1})\|
 &\ge c_k\rho_k\|\nabla_xL(\bar x,y_k)\|-C\rho_k^2 \ge(c_k-2C^2)\rho_k\|\nabla_xL(\bar x,y_k)\|.
\end{aligned}
$$
Choose $c_*>\max\{1,2C^2+4C^2/\varepsilon\}$. For $c_k\in [c_*, (2\rho_k)^{-1})$, these bounds imply
$$
 \|x_{k+1}-\bar x\|
 \le\frac{2C^2}{c_k-2C^2}\|\nabla_xL(\bar x,y_{k+1})\|
 \le\frac{\varepsilon}{2}\|\nabla_xL(\bar x,y_{k+1})\|.
$$
This proves~\eqref{relax:cone-invariance}. The proof is complete.
\end{proof}

Under the conditions of Lemmas~\ref{relax:mixed-lower} and~\ref{relax:corrected-step},  we have
$$
 \rho_{k+1}
 \le C_7(1+c_k)\rho_k^2
 \le C_6C_7(1+c_k)e_{k+1}.
$$
For sufficiently large $c_k$, condition~\eqref{eq:rel} also holds at the next iterate, so the same comparison can be applied again. These estimates will be used to show that the doubling rule eventually stops. The following theorem first establishes local convergence using $\omega_k\le1/2$, then proves that $c_k$ is eventually constant, and finally derives the quadratic rates.

\begin{theorem}\label{relax:adaptive}
Suppose that SOSC holds at the reference KKT pair $\bar z$ and that Assumption~\ref{ass:lipschitz} holds. Fix $b_0>0$ and a sufficiently large $\ell_0>0$. For $z_0$ sufficiently close to $\bar z$ with $\mu_0\ge0$, Algorithm~\ref{alg:corrected-ssqp} is well defined and doubles $b_k$ only finitely many times. It either terminates at a KKT pair or generates a sequence $\{z_k\}$ converging to $z_\infty:=(\bar x,y_\infty)$, with $y_\infty\in\mathcal M_*$. In the latter case, there exists a constant $C_8>0$ such that
\begin{equation}\label{relax:rates}
 e_{k+1}\le C_8e_k^2,\qquad r_{k+1}\le C_8r_k^2
\end{equation}
for all sufficiently large $k$, where $r_k:=\|z_k-z_\infty\|$.
\end{theorem}

\begin{proof} \ 
Under SOSC, Assumption~\ref{ass:eb} holds by~\cite[Proposition~3 and Remark~2]{izmailov2012stabilized}. Fix $C\ge\max\{1,C_5,C_7\}$, where $C_5$ and $C_7$ are from Lemmas~\ref{relax:mixed-predictor} and~\ref{relax:corrected-step}, respectively, and let $M$ be as in~\eqref{eq:local-smoothness-bounds}.  If $\rho_k=0$, then $z_k$ is a KKT pair and $e_k\le\kappa\rho_k=0$. We therefore consider a sequence that does not terminate.

\medskip\noindent {\bf Step 1: Local convergence.} Proposition~\ref{prop:local-solvability}{\rm (ii)} ensures the existence of $\widehat{z}_k$ satisfying~\eqref{relax:predictor}--\eqref{relax:localization} for $z_k$ sufficiently close to $\bar{z}$. By Lemma~\ref{relax:mixed-predictor}{\rm (i)}, $\mathcal I_k=\mathcal A$ and $\widehat\mu_{k,\mathcal N}=0$. From~\eqref{relax:update} and~\eqref{relax:localization}, we obtain
$$
 \mu_{k+1}\ge0,\qquad \mu_{k+1,\mathcal N}=0,\qquad \|z_{k+1}-z_k\|
 \le\|\widehat z_k-z_k\|
 \le\ell_0\rho_k.
$$

For $k\ge1$, we have $\mu_{k,\mathcal N}=0$. Set $E_k:=e_k+\|\nabla_xL(\bar x,y_k)\|.$ Lemma~\ref{relax:mixed-predictor}{\rm (iii)} gives $C^{-1}\rho_k\le E_k\le C\rho_k$. By~\eqref{relax:eta-recurrence}, Lemma~\ref{relax:mixed-predictor}{\rm (ii)}, and $\omega_k\le1/2$,
\begin{equation}\label{relax:adaptive-bootstrap}
 E_{k+1} \le\tfrac12E_k+C\rho_k^2 \le\bigl(\tfrac12+C^2\rho_k\bigr)E_k \le\tfrac34E_k,
\end{equation}
where the neighborhood is chosen so that $C^2\rho_k\le1/4$. Consequently,
\begin{equation}\label{eq:rhok}
 \rho_j\le C^2(3/4)^{j-k}\rho_k, \qquad j\ge k\ge1,
\end{equation}
as long as the iterates remain in this neighborhood.

Choose $\delta>0$ so that the above local estimates hold on $\mathbb B(\bar z,\delta)$. For $z_0$ sufficiently close to $\bar z$, the first step is well defined and $\|z_1-\bar z\|+4\ell_0C^2\rho_1<\delta$. The step bound and~\eqref{eq:rhok} then give, by induction,
$$
\begin{aligned}
 \|z_k-\bar z\|+\|\widehat z_k-z_k\| \le\|z_1-\bar z\|+\ell_0\sum_{j=1}^{k}\rho_j \le\|z_1-\bar z\|+4\ell_0C^2\rho_1<\delta,
 \qquad k\ge1.
\end{aligned}
$$
Thus $z_k$ and $\widehat z_k$ remain in this neighborhood. By~\eqref{eq:rhok}, $\sum_{k=1}^{\infty}\rho_k<\infty$. Since $\|z_{k+1}-z_k\|\le\ell_0\rho_k$, the sequence $\{z_k\}$ is Cauchy. Let $z_\infty$ denote its limit. By~\eqref{relax:adaptive-bootstrap}, $x_k\to\bar x$, and $\rho_k\to0$ gives $\mathcal R(z_\infty)=0$ by continuity. Hence $z_\infty=(\bar x,y_\infty)$ with $y_\infty\in\mathcal M_*$. Moreover, the step bound \eqref{relax:localization} and~\eqref{eq:local-smoothness-bounds} imply
\begin{equation}\label{relax:full-tail}
\begin{aligned}
 r_k&\le\ell_0\sum_{j\ge k}\rho_j\le4\ell_0C^2\rho_k, \qquad \rho_k=\|\mathcal R(z_k)-\mathcal R(z_\infty)\|\le Mr_k.
\end{aligned}
\end{equation}

\medskip\noindent \textbf{Step 2. Finitely many doublings of $b_k$.} For $i\ge1$, since $x_{i+1}=\widehat x_i$, Lemma~\ref{relax:mixed-predictor}{\rm (ii)} gives $e_{i+1}\le C_5\rho_i^2$. Thus the triangle inequality yields
\begin{equation}\label{relax:step-comparison}
 C_5\rho_i^2\le\tfrac12e_i
 \quad\Longrightarrow\quad
 s_i\ge e_i-e_{i+1}\ge\tfrac12e_i.
\end{equation}

Fix $\varepsilon\in(0,1]$ as in Lemma~\ref{relax:mixed-lower}, set $\alpha:=C_6^{-1}$, and let $c_*$ be the corresponding threshold in Lemma~\ref{relax:corrected-step}. Choose
$$
 \bar c:=\max\{c_*,1,\sqrt{2C(1+\varepsilon^{-1})},4C/\alpha\},
 \qquad c_{\max}:=4\max\{\bar c,b_0\}.
$$
By~\eqref{eq:rhok}, $\sup_{i\ge1}\rho_i\le C^2\rho_1$. We may therefore take $z_0$ sufficiently close to $\bar z$ that
\begin{equation}\label{eq:rhoksmall}
 c_{\max}\rho_i<\tfrac12,\qquad
 C^2\max\{1+\varepsilon^{-1},(1+c_{\max})/\alpha\}\rho_i
 \le\tfrac12
 \qquad(i\ge1).
\end{equation}

If $b_k<\bar c$ for every $k\ge1$, only finitely many doublings are possible. Otherwise, let $j\ge1$ be the first index with $b_j\ge\bar c$. The update rule gives
$$
 \bar c\le b_j\le 2\max\{\bar c,b_0\}=c_{\max}/2.
$$
We show that at most one further doubling can occur.

Whenever $k\ge j$ and~\eqref{eq:rel} fails, Lemma~\ref{relax:mixed-predictor}{\rm (iii)} gives
$$
 \rho_k\le C_5\bigl(e_k+\|\nabla_xL(\bar x,y_k)\|\bigr)
 \le C(1+\varepsilon^{-1})e_k.
$$
Since $C_5\le C$, this bound and~\eqref{eq:rhoksmall} give $C_5\rho_k^2\le C^2(1+\varepsilon^{-1})\rho_k e_k\le e_k/2$. Applying~\eqref{relax:step-comparison} at $i=k$, we obtain
$$
 \rho_k\le 2C(1+\varepsilon^{-1})s_k
 \le \bar c^2s_k\le b_k^2s_k,
$$
so $b_{k+1}=b_k$. If~\eqref{eq:rel} fails at every $k\ge j$, the proof of this step is complete.

Now let $k\ge j$ be the first index satisfying~\eqref{eq:rel}. Then $b_k=b_j$. A doubling may occur at this index, but $c_k\in\{b_j,2b_j\}$ gives $\bar c\le c_k\le c_{\max}$. Hence $\omega_k=c_k\rho_k<1/2$ by~\eqref{eq:rhoksmall}. By using Lemmas~\ref{relax:mixed-lower} and~\ref{relax:corrected-step}, we have
$$
 e_{k+1}\ge\alpha\rho_k^2,\qquad
 \rho_{k+1}\le C_7(1+c_k)\rho_k^2
 \le \frac{C}{\alpha}(1+c_k)e_{k+1}.
$$
Moreover, \eqref{relax:cone-invariance} ensures that \eqref{eq:rel} holds at $k+1$.

Since $c_k\le c_{\max}$, the bound on $\rho_{k+1}$ above and~\eqref{eq:rhoksmall} imply $C_5\rho_{k+1}^2\le (C^2/\alpha)(1+c_k)\rho_{k+1}e_{k+1} \le e_{k+1}/2$. Applying~\eqref{relax:step-comparison} at $i=k+1$, we obtain
$$
 \rho_{k+1}\le \frac{2C}{\alpha}(1+c_k)s_{k+1}
 \le c_k^2s_{k+1}=b_{k+1}^2s_{k+1},
$$
where the second inequality uses $c_k\ge\max\{1,4C/\alpha\}$. Thus $c_{k+1}=b_{k+1}=c_k$. Since~\eqref{eq:rel} holds at $k+1$ and $c_{k+1}\in[\bar c,c_{\max}]$, the same argument applies at $k+1$. By induction, no further doubling occurs. It follows that $b_k$ doubles only finitely many times and $c_k\le c_{\max}$ for all $k$.

\medskip\noindent {\bf Step 3: Quadratic convergence rates.} For all sufficiently large $k$, there is a fixed $b_\infty\in(0,c_{\max}]$ such that $b_k=c_k=b_\infty$ and $\rho_k\le b_\infty^2s_k$. Also, $\omega_k=b_\infty\rho_k<1/2$, so~\eqref{relax:accepted-estimates} and~\eqref{relax:full-tail} imply
$$
\begin{aligned}
 r_{k+1}&\le4\ell_0C^2\rho_{k+1}
 \le4\ell_0C^2C_7(1+b_\infty)\rho_k^2\le4\ell_0C^3(1+c_{\max})M^2r_k^2.
\end{aligned}
$$
For the primal rate, Lemma~\ref{relax:mixed-predictor}{\rm (ii)} and $s_k\le e_k+e_{k+1}$ yield
$$
\begin{aligned}
 \rho_k&\le b_\infty^2(e_k+e_{k+1}) \le b_\infty^2e_k+C_5b_\infty^2\rho_k^2
 \le b_\infty^2e_k+Cb_\infty^2\rho_k^2.
\end{aligned}
$$
For large $k$, $Cb_\infty^2\rho_k\le1/2$, so moving the last term to the left gives $\rho_k\le2b_\infty^2e_k$. Therefore,
$$
 e_{k+1}\le C_5\rho_k^2\le4Cb_\infty^4e_k^2.
$$
Taking $C_8=\max\{4Cc_{\max}^4,\,4\ell_0C^3(1+c_{\max})M^2\}$ proves~\eqref{relax:rates}.
\end{proof}

The proof uses SOSC to ensure Assumption~\ref{ass:eb} and the existence of a predictor satisfying~\eqref{relax:localization}. For problems only involving equality constraints, both properties also follow from noncriticality, so SOSC can be weakened as follows.

\begin{corollary}\label{relax:equality-noncritical}
Suppose that problem~\eqref{eq:nlp} has only equality constraints, Assumption~\ref{ass:lipschitz} holds, and $\bar\lambda$ is a noncritical multiplier associated with $\bar x$. Then all conclusions of Theorem~\ref{relax:adaptive} remain valid.
\end{corollary}
\begin{proof}\ 
Noncriticality implies Assumption \ref{ass:eb}. Moreover, Proposition~\ref{prop:local-solvability}{\rm (ii)} ensures the existence of $\widehat{z}_k$ satisfying  \eqref{relax:predictor}-\eqref{relax:localization}.  The proof of Theorem~\ref{relax:adaptive} uses SOSC only to obtain these two properties and therefore applies with the inequality components omitted.
\end{proof}

\begin{remark}
From Theorem~\ref{relax:adaptive} or Corollary~\ref{relax:equality-noncritical}, we have $\rho_k=O(s_k)$ and $\omega_k=O(\rho_k)$. For all sufficiently large $k$, the corrected iterates satisfy the pSQP equations~\eqref{eq:psqp} with
$$
\begin{aligned}
 \omega_1^k
 &=A_k^\top(\widehat\lambda_k-\lambda_{k+1})
   +G_k^\top(\widehat\mu_k-\mu_{k+1}),\\
 \omega_2^k&=-\rho_k\Delta\widehat\lambda_k,\\
 (\omega_3^k)_{\mathcal A}
 &=-g_{\mathcal A}(x_k)-G_{k,\mathcal A}d_k,
 \qquad (\omega_3^k)_{\mathcal N}=0.
\end{aligned}
$$
We obtain
$$
 \|\omega_1^k\|+\|\omega_2^k\|+\|\omega_3^k\|
 =O(\rho_k^2)=O(s_k^2).
$$
Since $\|\chi_k\|\le\|\omega_3^k\|$ by~\eqref{eq:chiomega}, these estimates imply $\nu_k=o(s_k)$ and $\|\omega_1^k\|=o(s_k)$, as required in Theorem~\ref{thm:effective}. They also satisfy the perturbation conditions in~\cite[Theorems~3.3--3.4]{fernandez2010sharp}. Here, the perturbation estimates needed to apply these results rely on the local convergence and eventual constancy of $b_k$ established above.

\end{remark}

\section{Numerical experiments}\label{sec:numerics}
In this section, we compare ordinary sSQP with Algorithm~\ref{alg:corrected-ssqp}, referred to as corrected sSQP, on the counterexample from Section~\ref{sec:counterexample} and on selected problems from the DEGEN collection~\cite{izmailov2009examples}. Both methods use the exact Lagrangian Hessian and the same subproblem solver. For corrected sSQP, we set $b_0=1$.

\subsection{The counterexample revisited}
For problem~\eqref{eq:example}, both methods start from
$$
 x_0=-t_0,\qquad y_0=0,\qquad  \lambda_0=t_0^2(1-t_0a_0),\qquad t_0=10^{-2}.
$$
To choose $a_0$, we approximate the function $q_*$ in Lemma~\ref{lem:curve} using $\mathcal G$ from Appendix~\ref{app:curve}. Starting from the constant function $q_0(t)\equiv\sqrt2-1$, we set
$$
 q_1=\mathcal Gq_0,\qquad q_2=\mathcal Gq_1,
 \qquad a_0=q_2(t_0)\approx0.4140489291519012.
$$
Thus $a_0$ approximates $q_*(t_0)$, and the initial point approximates a point on the invariant curve.  We evaluate~\eqref{eq:exact-map} in high-precision arithmetic.

\begin{figure}[!ht]
\centering
\includegraphics[width=\textwidth]{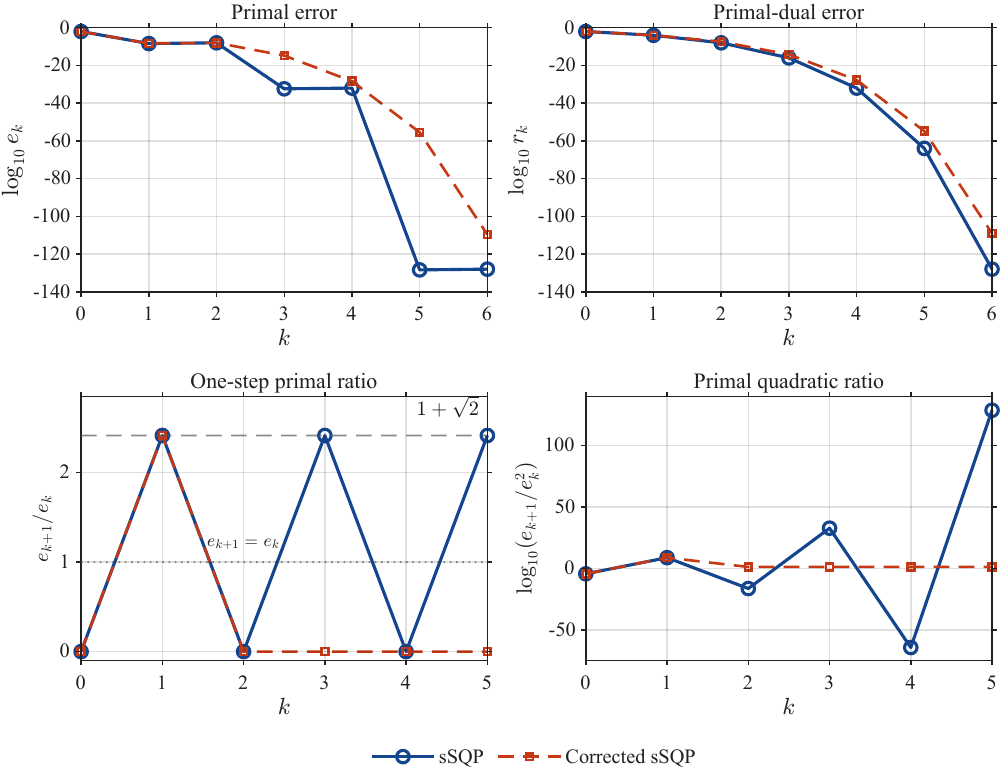}

\caption{Errors and convergence ratios for problem~\eqref{eq:example}.}
\label{fig:num-counterexample}
\end{figure}

Since the unique KKT pair is the origin, the primal and primal-dual errors in this experiment are
$$
 e_k=\sqrt{x_k^2+y_k^2},\qquad
 r_k=\sqrt{x_k^2+y_k^2+\lambda_k^2},
$$
respectively. The numerical results are recorded in Figure~\ref{fig:num-counterexample}.

Since the initial point only approximates a point on the invariant curve, Figure~\ref{fig:num-counterexample} reports the first six steps. We can see that both methods achieve high-precision primal-dual pair when $k=6$. For ordinary sSQP, $e_{k+1}/e_k$ is close to $1+\sqrt{2}$ at $k=1,3,5$, in agreement with Theorem~\ref{thm:counterexample}. For corrected sSQP, the primal error also increases at the second step, with $e_2/e_1\approx2.4138$. Indeed, after the first correction, $|\lambda_1|\approx10^{-4}$ and $e_1\approx4.14\times10^{-9}$, so~\eqref{eq:exact-map} gives $e_2\approx10^{-8}>e_1$. For $k\ge 3$, the ratios tends to zero, which agrees with the behavior proved in Theorem~\ref{relax:adaptive}. The lower-right panel of Figure~\ref{fig:num-counterexample} plots $\log_{10}(e_{k+1}/e_k^2)$ for both methods. For corrected sSQP, the last three values are close to $1.23$, while the corresponding values for ordinary sSQP oscillate with increasing magnitude.

\subsection{Comparison on DEGEN problems}
We next compare the two methods on 13 problems from the DEGEN collection, which originates from the examples of Izmailov and Solodov~\cite{izmailov2009examples}, and was also used in the numerical study of stabilized SQP by Izmailov, Solodov and Uskov~\cite[Section~4]{izmailov2015combining}. The dataset is available from the website\footnote{\url{https://w3.impa.br/~optim/DEGEN_collection.zip}}.  The selected problems have known solutions and belong to three classes: problems satisfying LICQ but not strict complementarity, problems satisfying MFCQ but not LICQ, and problems violating MFCQ. The identifiers of these three classes begin with 0, 1, and 2, respectively. The selected problems are listed in Table~\ref{tab:num-reference}.

\begin{table}[!ht]
\centering
\caption{Reference KKT pairs for the selected DEGEN problems. }
\label{tab:num-reference}
\small
\begin{tabular}{lrrrll}
\toprule
Problem & $n$ & $m$ & $p$ & $\bar x$ & $\bar y=(\bar\lambda,\bar\mu)$ \\
\midrule
00201 & 2 & 0 & 2 & $(0,0)$ & $(1,0)$ \\
00301 & 3 & 1 & 4 & $(0,0,1)$ & $(-2,0,4,0,0)$ \\
00302 & 3 & 0 & 6 & $(3,3,3)$ & $0\in\mathbb R^6$ \\
10201 & 2 & 0 & 2 & $(0,0)$ & $(1/4,3/4)$ \\
10202 & 2 & 0 & 2 & $(0,0)$ & $(1/8,1/16)$ \\
10203 & 2 & 0 & 2 & $(0,0)$ & $(1/2,1/8)$ \\
10204 & 2 & 0 & 2 & $(0,0)$ & $(1/2,1/4)$ \\
10205 & 2 & 0 & 3 & $(0,0)$ & $(1/8,1/16,0)$ \\
10207 & 2 & 0 & 2 & $(0,0)$ & $(7/8,1/8)$ \\
20204 & 2 & 2 & 0 & $(0,0)$ & $(0,0)$ \\
20211 & 2 & 0 & 2 & $(0,0)$ & $(1/4,1/4)$ \\
20303 & 3 & 2 & 0 & $(0,0,0)$ & $(0,0)$ \\
20304 & 3 & 2 & 0 & $(0,0,0)$ & $(0,0)$ \\
\bottomrule
\end{tabular}
\end{table}

For each problem, we generate 20 random initial primal-dual pairs near the reference pair listed in Table~\ref{tab:num-reference}. For each pair, we generate independent standard Gaussian vectors $\xi\in\mathbb R^n$ and $\eta\in\mathbb R^{m+p}$, and set
$$
 x_0=\bar x+\sigma\frac{\xi}{\|\xi\|},\qquad
 (\lambda_0,\widetilde\mu_0)=\bar y+\sigma\frac{\eta}{\|\eta\|},
 \qquad \mu_0=\max\{\widetilde\mu_0,0\}.
$$
where $\sigma=10^{-2}$. Both methods use the same 20 initial pairs for each problem.

Both methods use double precision and set $\ell_0=100$ in~\eqref{relax:localization}. Since each test problem has at most six inequality constraints, we solve the stabilized KKT linear system for each possible active set and choose the KKT point with the smallest primal-dual step. A run succeeds if
$$
 \rho_k=\|\mathcal R(z_k)\|\le10^{-12}
$$
within 40 iterations. If the solver cannot find a KKT point satisfying~\eqref{relax:localization}, we record the run as unsuccessful.

We record the number of iterations (Iter), the elapsed time in seconds (Time), the primal error (Perr), and the dual error (Derr). At the final iteration $k$, the two errors are
$$
 \mathrm{Perr}:=\|x_k-\bar x\|,\qquad
 \mathrm{Derr}:=\operatorname{dist}(y_k,\mathcal M_*).
$$

The numerical comparisons are summarized in Table~\ref{tab:num-iterations}. Both ordinary and corrected sSQP met the stopping criterion from all 260 initial primal-dual pairs across the 13 DEGEN problems. The average iteration counts are 3.10 for ordinary sSQP and 3.76 for corrected sSQP. In these tests, the correction slightly increases the iteration counts or running times, and the final primal and dual errors do not show a consistent advantage for either method. For corrected sSQP, the number of parameter doublings ranged from zero to five and averaged 2.53 per run.

\begin{table}[!ht]
\centering
\caption{Numerical comparison on DEGEN problems.}
\label{tab:num-iterations}
\scriptsize
\setlength{\tabcolsep}{3pt}
\begin{tabular}{lrrrrrrrr}
\toprule
& \multicolumn{4}{c}{Ordinary sSQP} & \multicolumn{4}{c}{Corrected sSQP} \\
\cmidrule(lr){2-5}\cmidrule(lr){6-9}
Problem &Iter & Time &Perr & Derr & Iter & Time & Perr & Derr \\
\midrule
00201 & 3.00 & 5.11e-04 & 7.52e-15 & 1.21e-14 & 3.30 & 7.30e-04 & 1.93e-13 & 5.56e-13 \\
00301 & 3.85 & 1.71e-03 & 5.59e-14 & 1.12e-12 & 4.05 & 1.85e-03 & 1.28e-14 & 2.64e-13 \\
00302 & 3.00 & 4.73e-03 & 5.65e-16 & 8.50e-17 & 3.70 & 5.78e-03 & 7.85e-16 & 1.09e-14 \\
10201 & 2.85 & 3.66e-04 & 8.03e-14 & 5.72e-17 & 3.15 & 3.89e-04 & 6.59e-14 & 8.53e-14 \\
10202 & 3.00 & 3.20e-04 & 2.29e-14 & 7.36e-17 & 4.75 & 5.05e-04 & 1.03e-17 & 1.49e-15 \\
10203 & 3.00 & 3.24e-04 & 1.27e-14 & 3.36e-16 & 3.90 & 4.26e-04 & 5.61e-15 & 7.17e-15 \\
10204 & 3.00 & 3.46e-04 & 2.50e-16 & 2.86e-16 & 3.75 & 4.17e-04 & 4.07e-16 & 7.05e-15 \\
10205 & 3.15 & 6.48e-04 & 7.28e-14 & 3.61e-17 & 4.55 & 9.23e-04 & 6.95e-17 & 6.39e-15 \\
10207 & 3.40 & 3.35e-04 & 5.81e-14 & 4.09e-17 & 3.45 & 3.44e-04 & 8.17e-14 & 6.93e-14 \\
20204 & 3.00 & 1.05e-04 & 2.74e-17 & 3.30e-17 & 3.35 & 1.34e-04 & 1.42e-14 & 9.01e-14 \\
20211 & 3.00 & 2.92e-04 & 4.20e-17 & 0.00e+0 & 3.10 & 3.05e-04 & 7.52e-14 & 1.35e-13 \\
20303 & 3.00 & 1.05e-04 & 8.34e-16 & 1.09e-15 & 3.90 & 1.37e-04 & 5.22e-18 & 5.60e-17 \\
20304 & 3.00 & 1.14e-04 & 1.07e-16 & 9.60e-17 & 3.90 & 1.42e-04 & 2.95e-17 & 5.98e-16 \\
\bottomrule
\end{tabular}
\end{table}

\bibliographystyle{siam}
\bibliography{references}

\appendix
\section{A uniform bound for the sSQP step}\label{app:small-step}
This appendix proves Proposition~\ref{prop:local-solvability}{\rm (i)} under Assumption~\ref{ass:eb}. We first use a projection property of the multiplier update to control the multiplier step by the primal step, and then use this assumption to estimate both steps.

\noindent
{\bf Proof of Proposition~\ref{prop:local-solvability}{\rm (i)}:}\ 
Write $y:=(\lambda,\mu)$, $c:=(h,g)$, $J:=c'$, and $\bar J:=J(\bar x)$. Set $d:=x^+-x$, $s:=\|d\|$, and $e:=\|x-\bar x\|$. Consider any pair $z,z^+$ sufficiently close to $\bar z$ satisfying \eqref{eq:ssqp-system} with $\mu\ge0$ and $\rho>0$. By Assumption~\ref{ass:eb}, $e\le\kappa\rho$.

\medskip\noindent \textbf{Step 1:} For $i\in\mathcal N$, the stabilized linearized constraint is strictly negative when $z$ and $z^+$ are sufficiently close to $\bar z$, because $g_i(\bar x)<0$, so complementarity gives $\mu_i^+=0$. At the input, $\min\{\mu_i,-g_i(x)\}=\mu_i$ for such $i$, and hence $\|\mu_{\mathcal N}\|\le\rho$. Thus $y^+$ belongs to the polyhedron
$$
 \mathcal F:=\{v=(v_\lambda,v_\mu):
 \bar J^\top v=\bar J^\top y^+,\ 
 v_{\mu,\mathcal A}\ge0,\ v_{\mu,\mathcal N}=0\}.
$$

By defining  $ r:=c(x)+J(x)d-c(\bar x)-\bar J(x^+-\bar x),$ we have from the definition fo $x^+$ that
$$
 r=c(x)-c(\bar x)-\bar J(x-\bar x)+[J(x)-\bar J]d.
$$
Thus, it follows from the Taylor expansion and the local Lipschitz continuity of $J$ that $\|r\| \leq O(e^2 + es)$, which along with Assumption \ref{ass:eb} implies that $\|r\|/\rho=O(\rho+s)$. For every $v\in\mathcal F$, the constraint equations in~\eqref{eq:ssqp-system} give
$$
 \langle c(x)+J(x)d-\rho(y^+-y),v-y^+\rangle\le0.
$$
Here $\bar J^\top(v-y^+)=0$ by definition of $\mathcal F$, and $\langle c(\bar x),v-y^+\rangle=0$ since the active constraints vanish at $\bar x$ and both inactive multiplier blocks are zero. Substituting $r$ and dividing by $\rho$ yields $\langle y+r/\rho-y^+,v-y^+\rangle\le0$ for all $v\in\mathcal F$. The projection characterization therefore gives
\begin{equation}\label{eq:app-projection}
 y^+=P_{\mathcal F}(y+r/\rho).
\end{equation}

To estimate $\operatorname{dist}(y,\mathcal F)$, the stationarity equation in~\eqref{eq:ssqp-system} gives
$$
 \bar J^\top(y^+-y)
 =-\nabla_xL(z)-\nabla^2_{xx}L(z)d
   +[\bar J-J(x)]^\top(y^+-y).
$$
The three terms have norms $O(\rho)$, $O(s)$, and $O(e)$, respectively, since $\|\nabla_xL(z)\|\le\rho$ and $y^+-y$ is locally bounded. As $e\le\kappa\rho$ and $\mu\ge0$, Hoffman's error bound~\cite{hoffman2003approximate} yields
$$
 \operatorname{dist}(y,\mathcal F)
 \le K_H\bigl(\|\bar J^\top(y-y^+)\|+\|\mu_{\mathcal N}\|\bigr)
 =O(\rho+s).
$$
Finally, by~\eqref{eq:app-projection} and nonexpansiveness of $P_{\mathcal F}$,
\begin{equation}\label{eq:app-dual-control}
 \|y^+-y\|
 \le \|r\|/\rho+\operatorname{dist}(y,\mathcal F)
 \le C(\rho+s)
\end{equation}
for some $C>0$ independent of the pair $z,z^+$.

\medskip\noindent \textbf{Step 2.} We show that $s\le 4\kappa\rho$ for every sufficiently close pair satisfying~\eqref{eq:ssqp-system} with $\mu\ge0$ and $\rho>0$. Otherwise, there exist pairs $(z_j,z_j^+)$ satisfying these conditions such that $z_j,z_j^+\to\bar z$ and $s_j>4\kappa\rho_j$ for every $j$, where $s_j:=\|x_j^+-x_j\|$ and $\rho_j:=\|\mathcal R(z_j)\|$. Then $s_j\to0$, $\rho_j=O(s_j)$, and~\eqref{eq:app-dual-control} gives $\tau_j:=\|y_j^+-y_j\|=O(s_j)$.

In the following estimate, we omit the subscript $j.$ By the stationarity equation in~\eqref{eq:ssqp-system} and continuity of $\nabla^2_{xx}L$, we have
\begin{equation*}
\begin{aligned}
 \nabla_xL(z^+)
 & ={}\nabla_xL(x^+,y)-\nabla_xL(x,y) -\nabla_{xx}^2L(z)(x^+-x)+[J(x^+)-J(x)]^\top(y^+-y) \\
 & = [J(x^+)-J(x)]^\top(y^+-y) + o(s).
\end{aligned}
\end{equation*}
Taking norms and using local Lipschitz continuity of $J$, we obtain $\|\nabla_xL(z^+)\|\le O(s\tau)+o(s)$.

For the constraint residuals, Taylor expansion of $h$ and $g$, the remaining equations in~\eqref{eq:ssqp-system}, and nonexpansiveness of the componentwise minimum in its second argument yield
$$
 \|h(x^+)\|+\|\min\{\mu^+,-g(x^+)\}\|
 =O(s^2+\rho\tau).
$$
Since $\rho_j=O(s_j)$ and $\tau_j=O(s_j)$, it follows that $\|\mathcal R(z_j^+)\|=o(s_j)$.

Applying Assumption~\ref{ass:eb} at both $z_j$ and $z_j^+$, we obtain
$$
 s_j\le \|x_j-\bar x\|+\|x_j^+-\bar x\|
 \le \kappa\bigl(\rho_j+\|\mathcal R(z_j^+)\|\bigr)
 \le \tfrac14 s_j+o(s_j).
$$
This contradicts $s_j>0$. Thus $s\le4\kappa\rho$ throughout a sufficiently small neighborhood. By~\eqref{eq:app-dual-control},
$$
 \|z^+-z\|\le s+\|y^+-y\|  \le \bigl(4\kappa+C(1+4\kappa)\bigr)\rho,
$$
which proves~\eqref{eq:zbound}.

\section{The two-step return and its invariant curve}\label{app:curve}
This appendix proves the two lemmas in Section~\ref{sec:counterexample}.

\medskip
\begin{proof}[Proof of Lemma~\ref{lem:two-step}]
Set $(x,\lambda)=z(t,a)=(-t,t^2(1-ta))$. We write $O_{\rm sm}(t^m)$ for a term $t^m r(t,a)$, where $r$ is jointly smooth on a neighborhood of $[0,\delta_0]\times I$ and bounded on this rectangle. We take $\delta_0>0$ sufficiently small below.

For the first step, set $S=\sqrt{1+t^2(2-ta)^2}=1+O_{\rm sm}(t^2)$. By~\eqref{eq:xlambda}, $\rho(x,\lambda)=tS$ and $x^2-\lambda=at^3$. Substituting them into~\eqref{eq:exact-map}, we have
$$
 x^+=t^4\frac{aS}{1-2t^2S}
     =t^4\bigl(a+O_{\rm sm}(t^2)\bigr),\qquad
 \lambda^+=t^2+2tx^+
     =t^2\bigl(1+O_{\rm sm}(t^3)\bigr).
$$
This proves~\eqref{eq:first-step}. For the second step, using the preceding expansions,
$$
 \rho^+:=\sqrt{\bigl((x^+)^2+\lambda^+\bigr)^2+(x^+)^2}
       =t^2\bigl(1+O_{\rm sm}(t^3)\bigr),\qquad
 (x^+)^2-\lambda^+=-t^2\bigl(1+O_{\rm sm}(t^3)\bigr).
$$
Another application of~\eqref{eq:exact-map} yields
$$
 T(t,a):=-x^{++}  =-\frac{\rho^+\bigl((x^+)^2-\lambda^+\bigr)}{1+2x^+\rho^+} =t^4\bigl(1+t^3\theta(t,a)\bigr)
$$
for a smooth function $\theta$. This proves~\eqref{eq:T-expansion}.

It remains to write $\lambda^{++}$ in the form $T^2(1+R)$. Since $x^{++}=-T$, we have $\lambda^{++}=(x^+)^2+2x^+T$. Hence, for $t>0$,
$$
 R(t,a):=\frac{\lambda^{++}}{T^2}-1 =\left(\frac{x^+}{T}\right)^2+2\frac{x^+}{T}-1.
$$
The expressions of $x^+$ and $T$ yield
$$
 \frac{x^+}{T} =\frac{a+O_{\rm sm}(t^2)}{1+t^3\theta(t,a)} =a+O_{\rm sm}(t^2).
$$
Thus $R$ extends smoothly to $t=0$ and satisfies $R(t,a)=a^2+2a-1+t^2\eta(t,a)$ for a smooth function $\eta$. This proves~\eqref{eq:two-arrows} and~\eqref{eq:R-expansion}.

Finally, differentiating $T=t^4+t^7\theta$ gives $\partial_aT=O(t^7)$ and $\partial_tT=4t^3+O(t^6)$ uniformly for $a\in I$, since $\theta$ and its first derivatives are bounded on $[0,\delta_0]\times I$.
\end{proof}

\begin{proof}[Proof of Lemma~\ref{lem:curve}]
Choose a closed interval $J\subset I$ containing $c$ in its interior. By~\eqref{eq:R-zero}, we may take $J$ and $\delta>0$ sufficiently small that $\partial_aR\ge \sqrt{2}$ on $[0,\delta]\times J$ and $R(0,\cdot)$ has opposite signs at the endpoints of $J$. By~\eqref{eq:T-expansion}, we also have $0<T(t,a)\le 2t^4<t$ for $0<t\le\delta$ and $a\in J$.

Let $q\in C^1([0,\delta])$ take values in $J$, with $q(0)=c$ and $\|q'\|_\infty\le1$. Define
$$
 H_q(t,a):=R(t,a)+T(t,a)q(T(t,a)).
$$
Let $T:=T(t,a)$. We have
$$
 \partial_t H_q
 =\partial_t R+
   \bigl(q(T)+Tq'(T)\bigr)\partial_t T, \qquad \partial_a H_q
 =\partial_a R+
   \bigl(q(T)+Tq'(T)\bigr)\partial_a T.
$$
It is clear that $|q(T)+Tq'(T)| \leq \max_{b\in J}|b|+\delta$. Since $\partial_tR=O(t)$, $\partial_tT=O(t^3)$, and $\partial_aT=O(t^7)$ by Lemma~\ref{lem:two-step}, we may reduce $\delta$, independently of $q$, so that $\partial_aH_q\ge \sqrt{2}/2$ and $|\partial_tH_q|\le \sqrt{2}/2$. Moreover, $H_q(t,a)=R(0,a)+O(t^2)$ uniformly in $a$ and $q$. Thus $H_q(t,\cdot)$ has opposite signs at the endpoints of $J$ and has a unique root $(\mathcal Gq)(t)\in J$ for each $t\in[0,\delta]$. At $t=0$ the root is $c$. The implicit function theorem gives $\mathcal Gq\in C^1([0,\delta])$ and
$$
 |(\mathcal Gq)'(t)|
 =\left|\frac{\partial_tH_q(t,(\mathcal Gq)(t))}
               {\partial_aH_q(t,(\mathcal Gq)(t))}\right|
 \le1.
$$

Starting from $q_0\equiv c$, set $q_{n+1}=\mathcal Gq_n$. For $n\ge1$, fix $t\in[0,\delta]$ and set $a=q_n(t)$. By construction, $H_{q_{n-1}}(t,a)=H_{q_n}(t,q_{n+1}(t))=0$. Applying the mean value theorem to $H_{q_n}(t,\cdot)$ between $a$ and $q_{n+1}(t)$ and using $\partial_aH_{q_n}\ge\sqrt{2}/2$, we obtain
$$
\begin{aligned}
 \frac{\sqrt{2}}{2}|q_{n+1}(t)-q_n(t)|
 \le |H_{q_n}(t,a)-H_{q_n}(t,q_{n+1}(t))| = |H_{q_n}(t,a)-H_{q_{n-1}}(t,a)|\le 2\delta^4\|q_n-q_{n-1}\|_\infty.
\end{aligned}
$$
Choose $\delta$ so that $\alpha:=4\delta^4/\sqrt{2}<1$. Taking the supremum over $t$ gives
$$\|q_{j+1}-q_j\|_\infty\le\alpha^j\|q_1-q_0\|_\infty, \quad {\rm for}\  j\ge0.$$
Hence $\{q_n\}$ is Cauchy in $C([0,\delta])$ and converges uniformly to some $q_*$. Since $\|q_n'\|_\infty\le1$, we have $|q_n(t)-q_n(s)|\le |t-s|$ for all $s,t\in[0,\delta]$. Letting $n\to\infty$ gives the same bound for $q_*$, so $\operatorname{Lip}(q_*)\le1$. Also, $q_*(0)=\lim_{n\to\infty}q_n(0)=c$.  Fix $t\in[0,\delta]$. By continuity of $T$, $T_n:=T(t,q_{n+1}(t))\to T_*:=T(t,q_*(t))$. Thus
$$
 |q_n(T_n)-q_*(T_*)|
 \le |q_n(T_n)-q_n(T_*)|+|q_n(T_*)-q_*(T_*)|
 \le |T_n-T_*|+\|q_n-q_*\|_\infty\to0.
$$
Since $R$ is continuous, passing to the limit in $R(t,q_{n+1}(t))+T_nq_n(T_n)=0$ gives $R(t,q_*(t))+T_*q_*(T_*)=0$, as required in~\eqref{eq:invariance}.

Finally, set $t_+=T(t,q_*(t))$. Then $0<t_+<t$ for $t>0$, and \eqref{eq:invariance} gives $R(t,q_*(t))=-t_+q_*(t_+)$. Substitution in~\eqref{eq:two-arrows} yields $\Phi^2(z(t,q_*(t)))=z(t_+,q_*(t_+))$, proving~{\rm (i)}. Part~{\rm (ii)} follows from~\eqref{eq:T-expansion}.
\end{proof}

\end{document}